\documentclass[a4paper,final,11pt]{amsart}

\usepackage{amsmath,amssymb,exscale,amsfonts,amssymb}
\usepackage[foot]{amsaddr}
\usepackage{times}
\usepackage{mathrsfs}
\usepackage{ams}
\usepackage{graphicx}
\usepackage{color, float, cases}
\usepackage{enumerate}
\usepackage{pgfplots}
\usepackage{tikz}
\usetikzlibrary{calc,fadings,decorations.pathreplacing,arrows,
           datavisualization.formats.functions,shapes.geometric}
\usepackage{verbatim}
\usepackage[margin=1in]{geometry}

\usepackage[utf8]{inputenc}
\usepackage{dutchcal}
\usepackage[displaymath, mathlines]{lineno}

\def\Dw{{\mathcal D}(\omega)}
\def\DD{{\mathcal D}_{D}(\omega)}
\def\DN{{\mathcal D}_{N}(\omega)}

\definecolor{darkolivegreen}{rgb}{0.33, 0.42, 0.18}
\definecolor{camouflagegreen}{rgb}{0.47, 0.53, 0.42}
\definecolor{calpolypomonagreen}{rgb}{0.12, 0.3, 0.17}
\definecolor{cadmiumgreen}{rgb}{0.0, 0.42, 0.24}
\definecolor{ao(english)}{rgb}{0.0, 0.5, 0.0}

\newcommand{\cor}[1]{\textcolor{black}{#1}}

\newcommand{\corg}[1]{\textcolor{black}{#1}}
\newcommand{\corj}[1]{\textcolor{black}{#1}}
\newcommand{\corb}[1]{\textcolor{black}{#1}}

\newcommand{\mcA}{{\mathcal A}}
\newcommand{\mcB}{{\mathcal B}}
\newcommand{\mcC}{{\mathcal C}}

\newcommand{\mcF}{\mathcal{F}}
\newcommand{\mcG}{\mathcal{G}}
\newcommand{\mcI}{{\mathcal I}}
\newcommand{\mcJ}{{\mathcal J}}

\newcommand{\mcP}{{\mathcal P}}
\newcommand{\mcQ}{{\mathcal Q}}
\newcommand{\mcR}{{\mathcal R}}
\newcommand{\mcS}{\mathcal{S}}
\newcommand{\mcT}{{\mathcal T}}

\newcommand{\mcZ}{\mathcal{Z}}

\newcommand{\msF}{\mathscr{F}}

\newcommand{\Real}{\mathop{\text{\rm Re}}}
\newcommand{\Imag}{\mathop{\text{\rm Im}}}

\newcommand{\jj}{\mathbf{j}}
\newcommand{\pp}{\mathbf{p}}

\newcommand{\bY}{\mathbf{Y}}
\newcommand{\bQ}{\mathbf{Q}}
\newcommand{\bb}{\mathbf{b}}

\newcommand{\bgamma}{\boldsymbol{\gamma}}

\newcommand{\bmeta}{\boldsymbol{\beta}}

\newcommand{\bq}{{\bf q}}
\newcommand{\bn}{{\bf n}}
\newcommand{\bfq}{{\bf f}}
\newcommand{\bg}{{\bf g}}

\newcommand{\bv}{{\bf v}}
\newcommand{\bd}{{\bf d}}

\newcommand{\ba}{{\bf a}}

\newcommand{\bw}{{\cor{\varkappa}}}

\newcommand{\bs}{{\bf s}}
\newcommand{\by}{{\bf y}}
\newcommand{\bx}{{\bf x}}
\newcommand{\bz}{{\bf z}}

\newcommand{\bC}{{\bf C}}
\newcommand{\bJ}{{\bf J}}
\newcommand{\bP}{{\bf P}}
\newcommand{\bG}{{\bf G}}
\newcommand{\bV}{{\bf V}}

\newcommand{\ii}{\mathbf{i}}

\newcommand{\mm}{\mathbf{m}}

\newcommand{\0}{\mathbf{0}}

\newcommand{\zzero}{\boldsymbol{0}}
\newcommand{\oone}{\boldsymbol{1}}

\newcommand{\bPhi}{\boldsymbol{\zeta}}

\newcommand{\Bset}{{\mathbb B}}
\newcommand{\Kset}{{\mathbb K}}

\newcommand{\Nset}{{\mathbb N}}

\newcommand{\Pol}{\mathbb{P}}
\newcommand{\eset}[1]{{\mathbb E} \left[ #1 \right] }

\newcommand{\var}{\text{\rm var}}

\newtheorem{assumption}{Assumption}

\newtheorem{remark}{Remark}
\newtheorem{problem}{Problem}

\newtheorem{theorem}{Theorem}
\newtheorem{lemma}{Lemma}

\DeclareMathOperator*{\esssup}{ess\,sup}
\DeclareMathOperator*{\essinf}{ess\,inf}

\begin{document}

\title[A stochastic collocation approach for parabolic PDEs with
    random domain deformations]
{A stochastic collocation approach for parabolic PDEs with
    random domain deformations}

\author{Julio E. Castrill\'{o}n-Cand\'{a}s}
\address{Boston University,
Department of Mathematics and Statistics,
111 Cummington Mall,
Boston, MA 02215
}
\email{jcandas@bu.edu}

\thanks{This material is based upon work supported by the National
  Science Foundation under Grant No. 1736392.
Research reported in this technical report was supported in part by the
National Institute of General Medical Sciences (NIGMS) of the
National Institutes of Health under award number 1R01GM131409-01.}

\author{Jie Xu}
\email{xujie@bu.edu}


\maketitle

\begin{abstract}
This work considers the problem of numerically approximating
statistical moments of a Quantity of Interest (QoI) that depends on
the solution of a linear parabolic partial differential equation.  The
geometry is assumed to be random and is parameterized by $N$ random
variables.  The parabolic problem is remapped to a fixed deterministic
domain with random coefficients and shown to admit an extension on a
well defined region embedded in the complex hyperplane. A Stochastic
collocation method with an isotropic Smolyak sparse grid is used to
compute the statistical moments of the QoI. In addition, convergence
rates for the stochastic moments are derived and compared to numerical
experiments.
\end{abstract}


\section{Introduction}

Mathematical modeling forms an essential part for understanding many
engineering and scientific applications with physical domains. These
models have been widely used to predict the QoI of any particular
problem when the underlying physical \cor{phenomenon} is well
understood.  However, in many cases the practicing engineer or
scientist does not have direct access to the underlying geometry and
uncertainty is introduced. It is essential to quantify the influence
of the domain uncertainty on the QoI.

In this paper a numerical method to efficiently solve parabolic PDEs
with respect to random geometrical deformations is developed.
Application examples include subsurface aquifers with soil variability
diffusion problems, ocean wave propagation (sonar) with geometric
uncertainty, chemical diffusion with uncertain geometries, among
others.

Collocation and perturbation approaches have been developed to
quantify the statistics of the QoI for elliptic PDEs with random
domains. The perturbation approaches
\cite{Harbrecht2008,Zhenhai2005,Guignard2017} are accurate for small
domain perturbations. In contrast, the collocation approaches
\cite{Chauviere2006,Fransos2008,Tartakovsky2006} allow the computation
of the statistics for larger domain deviations, but lack a full error
convergence analysis.  In \cite{Castrillon2016} the authors present a
collocation approach for elliptic PDEs based on Smolyak grids. An
analyticity analysis is \cor{performed. Convergence rates are derived
  and compared with numerical experiments.} \cor{Similar results where
  also obtained by the authors in \cite{Harbrecht2016,Hiptmair2018}}.

\cor{For stationary Stokes and Navier-Stokes Equations for viscous
  incompressible flow in \cite{Cohen2018}, a regularity analysis of the
  solution is studied with respect to the deformation of the
  domain. This approach is similar to the mapping technique proposed
  in this paper i.e. the random domain is assumed to be transformed
  from a fixed reference domain. The authors establish shape
  holomorphy with respect to the transformations of the shape of the
  domain.}

\cor{In \cite{Jerez-Hanckes2017} a shape holomorphy analysis for
  time-harmonic, \corg{electromagnetic} fields arising from scattering
  by perfect conductor and dielectric bounded obstacles. This approach
  \corg{falls} under the class of asymptotic methods for arbitrarily close
  random perturbations of the geometry. However, the authors show
  dimension-independent convergence rates for shape Taylor expansions
  of linear and higher order moments.}

\cor{A fictitious domain approach combined with \corg{Wiener}
  expansions was developed in \cite{Canuto2007}, where the elliptic
  PDE is solved in a fixed domain.  In \cite{Nouy2007,Nouy2008} the
  authors introduce a level set approach to the random domain problem.
  In \cite{Scarabosio2017} a multi-level Monte Carlo has been
  developed. This approach is well suited for low regularity of the
  solution with respect to the domain deformations.  Related work on
  Bayesian inference for diffusion problems and electrical impedance
  tomography on random domains is considered in \cite{Gantner2018,
    Hyvonen2017}.}

The work developed in this paper is a extension of the analysis and
error estimates derived in \cite{Castrillon2016} to the parabolic PDE
setting with Neumann and Dirichlet boundary conditions. \cor{Moreover,
  the stochastic domain deformation representation is extended to a
  larger class of geometrical perturbations.  This class of
  perturbations was originally introduced in
  \cite{Harbrecht2016,Guignard2017}.}  A rigorous convergence analysis
of the collocation approach based on isotropic Smolyak grids is
presented. This consists of an analysis of the regularity of the
solution with respect to the stochastic domain parameters. It is then
shown that the solution can be analytically extended to a well defined
region in $\C^{N}$ with respect to the domain random variables.
Error estimates are derived both in the ``energy norm'' as well as on
functionals of the solution (\cor{Quantities of Interest}) for Clenshaw
Curtis abscissas that can be easily generalized to a larger class of
sparse grids.

The outline of the paper is as follows: In Section \ref{setup} the
mathematical problem formulation is discussed.  The random domain
parabolic PDE problem is remapped onto a deterministic domain with
random matrix coefficients. In Section \ref{analyticity} the solution
of the parabolic PDE is shown to be analytically \cor{extendable} on a
well defined region in $\C^{N}$.  In Section
\ref{stochasticcollocation} the stochastic collocation method and
sparse grids are introduced.  In Section \ref{erroranalysis} error
estimates for the mean and variance of the QoI with respect to the
sparse grid and truncation approximations are derived.  Finally, in
section \ref{numericalresults} numerical examples are presented.

\section{Problem setting}
\label{setup}

Let $\Dw \subset \mathbb{R}^{d}$ be an open bounded domain with
Lipschitz boundary $\partial \Dw$ that is dependent upon a
random parameter $\omega \in \Omega$, where $(\Omega, \mcF,
\mathbb{P})$ is a complete probability space. Here $\Omega$ is the set
of outcomes, $\mcF$ is a $\sigma$-algebra of events and
$\mathbb{P}$ is a probability measure.

Suppose that the boundary $\partial \Dw$ is split into two disjoints
sections $\partial \DD$ and $\partial \DN$.  Consider the
following boundary value problem such that the following equations
hold almost surely:

\begin{equation}
\begin{split}
\partial_t u(\cdot,t, \omega)
  -\nabla \cdot ( a(\cdot, \omega) \nabla u(\cdot,t, \omega) ) 
&=
  f(\cdot,t,\omega)\,\,\,\,\,\,\,\,\,\mbox{in}\,\, \Dw \times (0,T) \\ 
u(\cdot,t,\omega) 
& = \corj{0} \hspace{7mm}\mbox{on
    $\partial \DD \times (0,T)$} \\
a(\cdot,\omega) \nabla u(\cdot,t,\omega) \cdot \cor{\bn(\cdot,\omega)} 
&=g_N(\cdot,\omega)\hspace{7mm}\mbox{on
    $\partial \DN \times (0,T)$} \\
u(\cdot,0,\omega) &=u_0(\cdot)\hspace{11mm}\mbox{on
$\Dw \times \{ t = 0 \}$}
\end{split}
\label{setup:strongformulation}
\end{equation}
\corg{where $T >0$. Let $\mcG := \cup_{\omega \in \Omega} \Dw$, then
  the functions $a:\mcG \rightarrow \R$, $f:\mcG \times (0,T)
  \rightarrow \R$, and $u_0:\mcG \rightarrow \R$ are defined over the
  region of all the stochastic perturbations of the domain $\Dw$ in
  $\R^{d}$. Similarly, let $\partial \mcG := \cup_{\omega \in \Omega}
  \partial \Dw \subset \R^{d}$, \corj{then the boundary conditions
  $g_N:\partial \mcG \rightarrow
  \R$ are defined over all the stochastic perturbations of the
  boundary $\partial \Dw$.}}

Before the weak formulation is posed, some notation and definitions
are established.  Define $L^{q}_{P}(\Omega)$, $q \in [1, \infty]$, as
the space of random variables such that
\[
\begin{split}
  L^{q}_{P}(\Omega)
  &:=
  \{v \,\,|\,\, \int_{\Omega}
|v(\omega)|^{q}\,\mbox{d}\mathbb{P}(\omega) < \infty \} \,\,\mbox{and} \\
L^{\infty}_{P}(\Omega)
&:=
\{v \,\,|\,\, \corb{\mathbb{P} -} \esssup_{ \omega \in \Omega}
|v(\omega)| < \infty \},
\end{split}
\]
\noindent where $v:\Omega \rightarrow \mathbb{R}$ is a strongly
measurable \cor{function}.  For $M$ valued vector functions $\bv:D
\rightarrow \R^{M}$, $D \subset \R^{d}$, \cor{$\bv :=
  [v_1,\dots,v_M]$}, $1 \leq q < \infty$, let
\[
\begin{split}
  [L^{q}(D)]^{M}
  &:=
  \{\bv \,\,|\,\, \int_{D}
\sum_{n=1}^{M} |v_n(\bx)|^{q}\,\mbox{d}\bx< \infty \} \,\,\mbox{and} \\
[L^{\infty}(D)]^{M}
&:=
\{\bv \,\,|\,\, \esssup_{\bx \in D, n = 1, \dots, M }
|v_n(\bx)| < \infty \}.
\end{split}
\]

\corg{Let
\[
V(\Dw) := \{v \in H^1(\Dw) \,\,|\,\, \mbox{$v = 0$ on $\partial \DD$}\},
\]
and denote by $V^{*}(\Dw)$ the dual space of $V(\Dw)$.}

Let $\bY:=[Y_1, \dots, Y_{N}]$ be a $N$ valued random vector
measurable in $(\Omega, \mcF, \mathbb{P})$ taking values on
$\Gamma:=\Gamma_{1} \times \dots \times \Gamma_{N} \subset
\mathbb{R}^{N}$ and $\mcB(\Gamma)$ be the Borel $\sigma-$algebra.
Define the induced measure $\mu_{\bY}$ on $(\Gamma,\mcB(\Gamma))$ as
$\mu_{\bY} : = \mathbb{P}(\bY^{-1}(A))$ for all $A \in
\mcB(\Gamma)$. {Assuming that the induced measure is absolutely}
continuous with respect to the Lebesgue measure defined on $\Gamma$,
then there exists a density function $\rho(\by): \Gamma \rightarrow
[0, +\infty)$ such that for any event $A \in \mcB(\Gamma)$
\[
\mathbb{P}(\bY \in A) := \mathbb{P}(\bY^{-1}(A)) = \int_{A} \rho( \by
)\,\mbox{d} \by.
\]
Now, for any measurable function \corg{$\bQ \in
  [L^{1}_{P}(\Gamma)]^{N}$} the expected value is defined as
\[
\corg{\mathbb{E}[\bQ] := \int_{\Gamma} \by \rho( \by )\, \mbox{d}
\by.}
\]
For $q \in \Nset_{+}$ define the following spaces
\[
\begin{split}
  L^q(\Gamma) &:= \{
  \cor{\mbox{$v(\by):
    \Gamma \rightarrow \R$ is strongly measurable}}
\, | \, \int_{\Gamma} v(\by)^q \rho(\by) \mbox{d}\by < \infty  \}
\,
\mbox{and} \\
\,
L^{\infty}(\Gamma) &:= \{ \cor{\mbox{$v(\by):
    \Gamma \rightarrow \R$ is strongly measurable}}
\, | \, \corb{\rho(\by) \mbox{d}\by - } \esssup_{\by \in \Gamma} |v(\by)|
< \infty  \}.
\end{split}
\]


We now pose the weak formulation of equation
\eqref{setup:strongformulation} \cor{(See Chapter 7 in
  \cite{Evans1998} and Chapter 7 in \cite{Knabner2003})}:
\begin{problem}
\corj{Given that $f(\bx,t,\omega) \in L^{2}(0,T;$ $L^{2}(\Dw))$,
$g_N (\bx,\omega) \in
  L^{2}(\DN)$ and $u_0 \in L^{2}(\mcG)$ find $u(\bx,t,\omega) \in$
$L^{2}(0,T; V(\Dw))$, with $\partial_t u \in
  L^{2}(0,T;V^{*}({\mathcal D}(\omega)))$, s.t.}
\begin{equation}
\begin{split}
\int_{\Dw} 
\partial_t u v + 
a(\bx,\omega) \nabla u \cdot \nabla v\,\,\emph{d}\bx
&= l(\omega; v),\,\,\,\,\,\,\,\,\,\,\,\,\,
\mbox{in $\Dw \times (0,T)$} \\
u(\bx,t,\omega) 
&=0\hspace{17mm}\mbox{on
     $\partial \DD \times (0,T)$} \\
u(\bx,0,\omega) &=u_0 \hspace{15mm}\mbox{on
$\Dw \times \{ t = 0 \}$},
\end{split}
\end{equation}
$\forall v \in V(\Dw)$ almost surely,
\corj{and}
\[
\corj{\begin{split}
l(\omega;v) &:=
\int_{\Dw} f(\bx,t,\omega)v\,\, \emph{d}\bx
+ \int_{\partial \DN}
g_N(\bx,\omega) v\,\, \corg{\emph{d}S(\bx)}. \\
\end{split}}
\]
  \corj{Recall}
  that the Neumann boundary condition
  $g_N(\bx,\omega) \in L^{2}(\partial \Dw)$ is defined over $\partial
  \mcG$. 
\cor{Problem \ref{setup:Prob1} has a unique solution if the following
  assumption is satisfied (See \cite{Evans1998}, \cite{Knabner2003},
  \cite{Lions1972}):}
\label{setup:Prob1}
\end{problem}

\corb{
  \begin{remark}
    In Problem \ref{setup:Prob1} we assume vanishing Dirichlet
    boundary conditions.  We also considered a nonzero Dirichlet
    condition e.g. $ u(\cdot, t, \omega) = g_D(\cdot,\omega) \;
    \mbox{on $\partial \DD \times (0,T)$} $. For this setup there are
    several compatibility conditions for $g_{D}$ that must be
    satisfied.  First, certain regularity assumptions of $g_{D}$ have
    to be made and furthermore, it should follow that $a(\cdot,\omega)
    \nabla g_{D}(\cdot, \omega) \cdot \cor{\bn(\cdot,\omega)}
    =g_N(\cdot,\omega)$ on $(\partial \DD \cap \partial \DN) \times
    (0,T)$.  Second, considering the weak solution, as in Problem
\ref{setup:Prob1}, the integration by parts leads to an extra term of
the form $\int_{\Dw} \nabla \cdot (a(\bx, \omega)\nabla g_{D}(\bx,
\omega))v\,\, \emph{d}\bx$. Thus this extra term should be considered
in the analytic regularity analysis and error bounds described in this
paper. This is beyond the current scope of our work, as it is already
very extensive.  For simplicity, we set the Dirichlet condition to the
trivial condition.
\end{remark}
}

\begin{assumption} There exist constants $a_{min}$ and $a_{max}$ such
  that
\[
0 < a_{min} \leq a(\bx, \omega) \leq a_{max} < \infty
\,\,\,\mbox{for a.e. $\bx
    \in \Dw$, $\omega \in \Omega$},
\]
\noindent where
\[
a_{min} := \essinf_{\bx \in \Dw, \omega \in \Omega } a(\bx, \omega)
\,\,\,\,\mbox{and}\,\,\,\, a_{max} := \esssup_{\bx \in \Dw,
  \omega \in \Omega} a(\bx, \omega).
\]
\label{setup:Assumption1}
\end{assumption}

\subsection{Reformulation on a reference domain}

To simplify the analysis of Problem \ref{setup:Prob1} we remap the
solution $u \in H^{1}(\Dw)$ onto a non-stochastic fixed domain. This
approach has been applied in
\cite{Fransos2008,Castrillon2016,Harbrecht2016,Hiptmair2018,Guignard2017}
and we can then take advantage of the extensive theoretical and
practical work of PDEs with stochastic diffusion coefficients.

Assume that given any $\omega \in \Omega$ the domain $\Dw$ can be
mapped to a reference domain $U \subset \mathbb{R}^{d}$ with Lipschitz
boundary through a random map \cor{$F(\omega):\overline U \rightarrow
  \overline{\Dw}$}, where $F$ is assumed to be a \cor{bijection}. The
map $\bmeta \mapsto \bx$, $\overline U \rightarrow \overline \Dw$, is
written as
\[
\cor{\bmeta \mapsto \bx = F(\bmeta,\omega),}
\]
\cor{where $\bmeta$ are the coordinates for the reference domain $U$.}
See the cartoon example in Figure \ref{setup:fig1}.
\begin{assumption} Given a one-to-one map $F(\bmeta,\omega):
  \overline U \rightarrow \overline{\Dw}$ there exist constants
  $\F_{min}$ and $\F_{max}$ such that
\[
0<\F_{min} \leq \sigma_{min} (\partial F(\omega)) 
\,\,\mbox{and}\,\, \sigma_{max} (\partial F(\omega)) \leq
\F_{max} < \infty
\label{setup:Assumption2}
\]
almost everywhere in $U$ and almost surely in $\Omega$. \cor{Denoted by
$\sigma_{min} (\partial F(\omega))$ (and $\sigma_{max} ($ $\partial
F(\omega)))$ the minimum (respectively maximum) singular value of the
Jacobian matrix $\partial F(\omega)$}.
\end{assumption}

\begin{remark}
  \cor{The previous assumption implies that the Jacobian $|\partial
    F(\bmeta, \omega)| \in L^{\infty}(U)$ almost surely.}
\end{remark}

\begin{figure}
\begin{center}
\begin{tikzpicture}[scale=1] 
    \foreach \x in {90,...,-90} { 
    \pgfmathsetmacro\elrad{20*max(sin(\x),.3)}
    \pgfmathsetmacro\ltint{.9*abs(\x-55)/180}
    \pgfmathsetmacro\rtint{.9*(1-abs(\x+55)/180)}
    \definecolor{currentcolor}{rgb}{0, \ltint, \rtint}
    \draw[color=currentcolor,fill=currentcolor] 
        (xyz polar cs:angle=\x,y radius=0.75,x radius=1.5) 
        ellipse (\elrad pt and 20pt);
    \definecolor{currentcolor}{rgb}{0, \ltint, \rtint}
    \draw[color=currentcolor,fill=currentcolor] 
        (xyz polar cs:angle=180-\x,radius=0.75,x radius=1.5) 
        ellipse (\elrad pt and 20pt);
    }
\coordinate (O) at (3.5,0);
\coordinate (P) at (4.8,0);
\draw[->, >=latex, gray, line width=4 pt] (O) -- (P);
\node[above=16pt,scale=1] at (4,-1.5) {$F(\omega)$};
\node[above=16pt,scale=1] at (-1.5,1) {$U$};
\end{tikzpicture}
\begin{tikzpicture}[scale=1] 
    \foreach \x in {90,89,...,0} { 
    \pgfmathsetmacro\elrad{20*max(sin(\x),.3)}
    \pgfmathsetmacro\ltint{.9*abs(\x-55)/180}
    \pgfmathsetmacro\rtint{.9*(1-abs(\x+55)/180)}
    \pgfmathsetmacro\ztint{2*(abs(\x)/180)}
    \definecolor{currentcolor}{rgb}{\ztint,\ltint,\rtint}
    \draw[color=currentcolor,fill=currentcolor] 
        (xyz polar cs:angle=\x,y radius=1.5,x radius=1.5) 
        ellipse (\elrad pt and 20pt);
    \definecolor{currentcolor}{rgb}{\ztint,\ltint,\rtint}
    \draw[color=currentcolor,fill=currentcolor] 
        (xyz polar cs:angle=180-\x,radius=1.5,x radius=1.5) 
        ellipse (\elrad pt and 20pt);
    }

    \foreach \x in {0,...,-90} { 
    \pgfmathsetmacro\elrad{20*max(sin(\x),.3)}
    \pgfmathsetmacro\ltint{.9*abs(\x-55)/180}
    \pgfmathsetmacro\rtint{.9*(1-abs(\x+55)/180)}
    \definecolor{currentcolor}{rgb}{0, \ltint, \rtint}
    \draw[color=currentcolor,fill=currentcolor] 
        (xyz polar cs:angle=\x,y radius=.75,x radius=1.5) 
        ellipse (\elrad pt and 20pt);
    \definecolor{currentcolor}{rgb}{0, \ltint, \rtint}
    \draw[color=currentcolor,fill=currentcolor] 
        (xyz polar cs:angle=180-\x,radius=.75,x radius=1.5) 
        ellipse (\elrad pt and 20pt);
    }
    \pgfmathsetmacro\elrad{cos(-135)}
    \pgfmathsetmacro\xrad{1.5cm-20pt*\elrad}
    \pgfmathsetmacro\yrad{.75cm-20pt*sin(-135)}

    \node[above=16pt,scale=1] at (-2.3,0.7) {${\mathcal D}(\omega)$};
\end{tikzpicture}
\end{center}
\caption{Cartoon example of stochastic domain realization from a
    reference domain. This figure is modified from the TikZ tex code
    from {\it Smooth map of manifolds and smooth spaces} by Andrew
    Stacey.}
\label{setup:fig1}
\end{figure}

From the Sobolev chain rule (see Theorem 3.35 in \cite{Adams1975} or
page 291 in \cite{Evans1998}) it follows that for any $v \in H^{1}(\Dw)$
\begin{equation}
\nabla_{\Dw} v = \partial F^{-T} \nabla (v \circ F),
\label{setup:equivalence}
\end{equation}
where $\nabla_{\Dw}$ refers to the gradient on the domain $\Dw$,
$\nabla$ is the gradient on the reference domain $U$, and $(v \circ F)
\in H^{1}(U)$. Let
\[
V := \{v \in H^1(U): \mbox{$v = 0$ on \cor{$\partial U_D$}}\},
\]
where $\partial U$ is the boundary of $U$, $\partial U_{D} \subset
\partial U$ is the range of $F^{-1}$ with respect to the boundary
$\DD$, $\partial U_{N} \subset \partial U$ is the range of $F^{-1}$
with respect to the boundary $\DN$ and $\partial U_{D} \cup \partial
U_{N} = \partial U$. \corg{Furthermore, denote by $V^{*}$ the dual
  space of $V$.}

We can now show that:
\begin{lemma}
  Under Assumptions \ref{setup:Assumption2}
the following pairs of spaces are isomorphic
\begin{enumerate}[i)]
\item $L^{2}(\Dw) \cong L^{2}(U)$.
\item $H^{1}(\Dw) \cong H^{1}(U)$.
\item $L^{2}(0,T; L^{2}(\Dw)) \cong L^{2}(0,T; L^{2}(U))$.
\item $L^{2}(0,T; H^{1}(\Dw)) \cong L^{2}(0,T; H^{1}(U))$.
\item $L^{2}(\partial \Dw) \cong L^{2}(\partial U)$.
\item $L^{2}(0,T;V^{*}(\Dw)) \cong L^{2}(0,T;V^{*})$.
\item $H^{1/2}(\partial \Dw) \cong H^{1/2}(\partial U)$. 
\end{enumerate}
\label{setup:lemma0}
\end{lemma}
\begin{proof}
  ~
\begin{description}
\item [$i) - iv)$] From the Sobolev chain rule it is not hard to prove.
\item [$v)$] Suppose we have a disjoint finite covering $\mcT$ of the
  boundary $\partial U$ such that for each $\tau \in \mcT$ there
  exists a Lipschitz bijective mapping $\xi_{\tau}:B^0_r \rightarrow
  \tau$ (c.f. trace theorem proof, p. 258 in \cite{Evans1998} for
  details and \cite{Sauter2011}), where $ B^0_r := \{ \bx \in B_r
  \,|\, x_d = 0\}$ and $B_r \subset \R^{d}$ is a ball of radius $r$.
  \corb{In the following proof the Lipschitz mappings $\xi_{\tau} $,
    $\tau \in \mcT$, are assumed to be differentiable.  From the
    Radamacher Theorem \cite{Federer1969} every Lipschitz function
    is differentiable almost everywhere. Therefore without loss of
    generality we can replace the Lipschitz mappings $\xi_{\tau}$,
    $\tau \in \mcT$, with an equivalent differentiable version except
    for sets of measure zero.}
  For simplicity we shall perform the following analysis with respect
  to a single open set $\tau$ and mapping $\xi_{\tau} : B^0_r
  \rightarrow \tau$. Let $\bJ_{\tau}:=\{\partial_{x_i}
  \corj{\xi_{\tau_j}} \}_{1 \leq i \leq d-1}^{1 \leq j \leq d}$, then
  for any $v \in L^{2}(\partial U)$
  \[
  \int_{\tau} v^2\,\mbox{d}S = 
  \int_{B^0_r} (v \circ \xi_{\tau})^2 \,
(det(\bJ_{\tau}^{T} \bJ_{\tau}))^{\frac{1}{2}}
  \, \mbox{d}\bx'. 
  \]
  Now, $K_{\tau} = F(\tau, \omega)$ covers a portion of the boundary
  of $\partial \Dw$, then
  \[
  \int_{K_{\tau}} v^2\,\mbox{d}S = \int_{B^{0}_r} (v \circ F \circ
  \xi_{\tau})^2 \,
(det(\bJ_{F \circ \tau}^{T} \bJ_{F \circ \tau}))^{\frac{1}{2}}
, \mbox{d}\bx', 
  \]
  where $\bJ_{F \circ \tau} = \partial F(\cdot, \omega)\bJ_{\tau}$. It
  is not hard to show that for any vector $\bs \in \R^{d-1}$, where
  $\|\bs\|_{l^2} = 1$,
  \[
  \begin{split}
  &\sigma_{min}(\partial F(\cdot,\omega)^{T} \partial F(\cdot,\omega))
  \sigma_{min}(\bJ_{\tau}^{T}\bJ_{\tau})
  \leq 
  \bs^T
  \bJ_{\tau}^{T} \partial F(\cdot,\omega)
  ^{T} \partial F(\cdot,\omega)
  \bJ_{\tau}
  \bs \\
  &\leq
  \sigma_{max}(\partial F(\cdot,\omega)^{T}\partial F(\cdot,\omega))
  \sigma_{max}(\bJ_{\tau}^{T}\bJ_{\tau}).
  \end{split}
  \]
  The result follows.
\item [$vi)$]
  \corg{Suppose that $\xi \in V(\Dw)^{*}$, then $\| \xi \|_{V(\Dw)^{*}}$
    is equal to
  \[ 
   \sup_{\begin{array}{c} v \in V(\Dw) \\ \|v\|_{V(\Dw)} \leq 1
      \end{array}}
  |\xi(v)|
  = \sup_{\begin{array}{c} v \circ F \in V \\ C\|v \circ F\|_{V} \leq \|v\|_{V(\Dw)
      } = 1
      \end{array}}
  |\xi(v \circ F)|.
  \]
  The positive constant $C > 0$ is due to the fact that $H^{1}(\Dw)
  \cong H^{1}(U)$.  Let $\hat w = C(v \circ F)$, 
  \corj{then}
  \[
  \begin{split}
    \| \xi \|_{V(\Dw)^{*}}
    & 
  \leq \sup_{\begin{array}{c} \hat w \in V \\ \| \hat w \|_{V} \leq 1
      \end{array}}
  C^{-1}|\xi(\hat w)|=  C^{-1}
  \| \xi \|_{V^{*}}, \forall C > 0.
  \end{split}
  \]
  The converse is similarly proven.
  }
  
\item [$vii)$] The result follows by using $ii)$, the Trace Theorem
  and inverse Trace Theorem (Theorems 2.21 and 2.22 in
  \cite{Steinbach2007}).
\end{description}
\end{proof}
Note that analogous lemmas are proved in
\cite{Castrillon2016,Harbrecht2016}.

\smallskip
In the rest of the paper the terms a.s. and a.e. will be dropped
unless emphasis or disambiguation is needed.

For any $v,s \in H^{1}(U)$
\[
\begin{split}
B(\omega; s,v)
&:= 
  \int_{U}
  (a \circ F)(\bmeta,\omega)
  \nabla s
^{T}
\cor{\partial F^{-1}(\bmeta,\omega)\partial F^{-T}(\bmeta,\omega)}
\nabla v\,
|\partial F(\bmeta,\omega)|\,
\mbox{d}\bmeta.
\end{split}
\]
With a change of variables the boundary value problem is
remapped. However, we first deal with the case where

\begin{problem} Given that $(f \circ F)(\bmeta,t,\omega)
\in L^{2}(0,T;L^{2}(U))$, 
$ \hat g_N :=g_N \circ F$, 
and $\hat g_N \in
L^{2}(\partial U_N)$ find $\hat{u}(\bmeta,t,\omega) \in L^{2}(0,T;
V)$, \cor{with $\partial_t u \in L^{2}(0,T;V^{*})$},
s.t.
\[
\begin{split}
\int_{U} 
v |\partial F(\bmeta,\omega)| \partial_t \hat{u}(\bmeta,t,\omega) 
\,\mbox{d}\bmeta
+ 
B(\omega; \hat{u}, v)
&= \hat{l}(\omega; v), \hspace{12mm}
\mbox{in $U \times (0,T)$} \\
\hat{u}(\bmeta,t,\omega) 
&=0, \hspace{20mm}\mbox{on
    $\partial U_D \times (0,T)$} \\
\hat{u}(\bmeta,0,\omega) &=(u_0 \circ F)(\bmeta,\omega)\hspace{3mm}\mbox{on
$U \times \{ t = 0 \}$}
\end{split}
\]
$\forall v \in V$ almost surely, where
\[
\corj{\begin{split}
\hat{l}(\omega;v) 
&:=
\int_{U} (f \circ F)(\bmeta,\omega) 
| \partial F(\bmeta) |v\,
\emph{d}\bmeta
\\
&+
\corg{\sum_{\tau \in \mcT}
\int_{B^{0}_{r}}
(g_N \circ F)(\bmeta \circ \xi_{\tau}
,\omega)
(v \circ \xi_{\tau})} \\
& \corg{det ( \bJ_{\tau}^{T} \partial F(\bmeta \circ \xi_{\tau},\omega)
  ^{T} \partial F(\bmeta \circ \xi_{\tau},\omega)
  \bJ_{\tau})^{\frac{1}{2}}
\,\emph{d}\bx',}
\end{split}}
\]
\label{setup:Prob2}
\end{problem}
\corj{where 
  $T_{U}:H^{1/2}(\partial U) \rightarrow H^{1}(U)$ is a linear bounded
  operator such that $\forall \hat{g} \in H^{1/2}(\partial U)$,
  $T_{U}\hat{g} \in H^{1}(U)$ satisfies $(T_{U}\hat{g})|_{\partial U} = \hat{g}$.}

The weak solution $u \in H^{1}(\Dw)$ for the non-zero Dirichlet
boundary value problem is simply obtained as \corj{$u(\bx,\omega) = (\hat{u}
\circ F^{-1})(\bx, \omega). $} 

Now we have to be a little careful. The existence theorems from
\cite{Evans1998}, Chapter 7, do not apply directly to Problem
\ref{setup:Prob2} due to the $|\partial F(\bmeta,\omega)| \partial_t
\hat{u}$ term. Although the existence proof in \cite{Evans1998} can be
modified to incorporate this extended term, we direct our attention
to Theorem 10.9 in \cite{Brezis2010} from J. Lions \cite{Lions1972}.

\corg{Let $H$ (with norm $\|\cdot\|_{H}$) and $W$ (with norm
  $\|\cdot\|_{W}$ )} be Hilbert spaces with the associated dual spaces
$H^{*}$ and $W^{*}$ respectively. It is assumed that $W \subset H$
with dense and continuous injection so that
\[
W \subset H \subset W^{*}.
\]
For a.e. $t \in [0,T]$ suppose the bilinear form $\corg{A[t;
    \zeta,v]}: W \times W \rightarrow \corg{\R}$ satisfies the
following properties:

\begin{enumerate}[i)]

\item For every $\corg{\zeta,v} \in W$ the function $t \mapsto \corg{A[t;
  \zeta,v]}$ is measurable,

\item For all $\corg{\zeta,v} \in W$ $\corg{|A[t; w,v]|} \leq M \corg{\|\zeta\|_{W}}
  \|v\|_{W}$ for a.e. $t
  \in [0,T]$

\item For all $v \in W$ $\cor{A[t; v,v]} \geq \alpha \corg{\|v\|^{2}_{W}
  - C\|v\|_{H}^{2}}$
  for a.e. $t \in [0,T]$.
  
  \end{enumerate}
where $\alpha > 0$, $M$ and $C$ are constants.

\begin{theorem}{(J. Lions)} 
  Given a bounded linear functional $\corg{\mcZ} \in L^{2}(0,T;W^{*})$ and
  $u_0 \in H$, there exists a unique function $\hat u$ satisfying
  $\hat u \in L^{2}(0,T;W) \cap C([0,T];H),$ $\partial_t \hat u \in
  L^{2}(0,T;W^{*})$
\[
\langle \partial_t \hat u, v \rangle
+ A[t; \hat u,v] = \langle \corg{\mcZ},
v \rangle
\]
for a.e. $t \in (0,T)$, $\forall v \in W$, and $\hat u(0) = u_0$.
\label{setup:theorem1}
\end{theorem}
\begin{proof} See \cite{Lions1972}.
  \end{proof}

We can now use Theorem \ref{setup:theorem1} to show that there exists
a unique solution to \corg{Problems} \ref{setup:Prob1} and
\ref{setup:Prob2}. Let $W=V(\Dw)$ and $H = L^{2}(\Dw)$ then from
Theorem \ref{setup:theorem1} there \corg{exists} a unique solution $u
\in L^{2}(0,T;V(\Dw)$ for Problem \ref{setup:Prob1} such that
$\partial_t u \in L^{2}(0,T;V^{*}$ $(\Dw))$.  From Lemma
\ref{setup:lemma0} there \corg{is} an isomorphic map between $\hat u$
and $u$. Since there is a unique solution for Problem
\ref{setup:Prob1}, we conclude there \corg{exists} a solution $\hat u
\in L^{2}(0,T;V)$ for Problem \ref{setup:Prob2} such that $\partial_t
\hat u \in L^{2}(0,T;V^{*})$.  The last step is to confirm that it is
unique solution. This is done by checking $\hat u = 0$ is \corg{the}
solution whenever $\hat l(\omega;\cdot) = (u_0 \circ
F)(\cdot,0,\omega) = 0$.


\subsection{Stochastic domain deformation map}
\label{setup:domainparametrization}
The next step is to build a parameterization of the map
$F(\bmeta,\omega)$ from a set of random variables $Y_1,\dots,Y_N$
\corg{with} \cor{probability density function $\rho(\by)$}.  One
objective is to build a parameterization such that a large class of
stochastic domain deformations are represented.  Following the same
approach as in \cite{Guignard2017,Harbrecht2016}, without loss of
generality we assume that the map $F(\bmeta,\omega)$ has the finite
noise model
\[
F(\bmeta, \omega) := \bmeta + \sum_{n = 1}^{N}
\sqrt{\mu_n}\bb_n(\bmeta)Y_n(\omega).
\]
From the Doob-Dynkin Lemma the solution $\hat u$ to Problem
\ref{setup:Prob2} will be a function of the random variables $Y_1,
\dots, Y_N$.

This is a very general representation of the random domain
deformation. For example, such representation may be achieved by a
truncation of a \cor{Karhunen-Lo\`{e}ve (KL)} expansion of vector
random fields \cite{Harbrecht2016}. In general, the KL eigenfunctions
$\bb_l(\bmeta) \in [L^{2}(U)]^{d}$, which presents a problem, as the
KL expansion of the random domain may lead to large spikes and thus
most likely Problem \ref{setup:Prob2} will be ill-posed.  However,
under stricter regularity assumptions of the covariance function the
eigenfunctions will have higher regularity (see
\cite{Frauenfelder2005} for details).  We thus make the following
assumptions:
\begin{assumption}~
\begin{enumerate}
\item $\bb_{1},\dots, \bb_{N} \in [W^{1,\infty}(U)]^{d}$.
\item $\| \bb_n \|_{[L^{\infty}(U)]^{d}} = 1$ for $n = 1,\dots N$.
\item $\mu_1, \dots, \mu_{N}$ are monotonically decreasing.
\end{enumerate}
\label{setup:Assumption4}
\end{assumption}

From the stochastic model formulated in Section \ref{setup} the
Jacobian \cor{matrix} $\partial F$ is written as
\begin{equation}
\partial F(\bmeta,\omega) = I + 
\sum_{n=1}^{N} \sqrt{\mu_{n}}
\partial \bb_{n}(\bmeta)
Y_{n}(\omega).
\label{analyticity:eqn1}
\end{equation}

\section{Analyticity of the boundary value problem}
\label{analyticity}

In this section we show that the solution to Problem \ref{setup:Prob2}
can be analytically extended on a region $\Theta_{\beta}$ in $\C^{N}$
with respect to stochastic domain $\by \in \Gamma$. The size of the
region $\Theta_{\beta}$ is \cor{related} to the regularity of the
solution with respect to $\Gamma$.  \cor{This provides us a path to
  estimate the convergence rates of the stochastic moments by using a
  sparse grid approximation.}  \cor{In particular, the larger the size
  of the region $\Theta_{\beta}$, the faster the convergence rate of
  the sparse grid approximation will be.}

\begin{remark}
\corg{To simplify the analysis assume that $\Gamma$ is bounded in
  $\R^{N}$. Without loss of generality it can also be assumed that
  $\Gamma = [-1,1]^{N}$. However, $\Gamma$ can be extended to the
  non-bounded case by following the approach described in
  \cite{Babuska2010}.}
\end{remark}

We formulate the region
$\Theta_{\beta}$ by making the following assumption:
\begin{assumption}
\begin{enumerate}
\item There exists $0 < \tilde{\delta} < 1$ such that $\sum_{n=1}^{N}
  \sqrt{\mu_{n}} \| \partial \bb_n(\bmeta) \|_{2} 
  \leq 1 - \tilde{\delta}$ for
  all $\bmeta \in U$.
\end{enumerate}
\label{analyticity:assumption1}
\end{assumption}
For any $0 < \beta < \tilde{\delta}$ define the region $\Theta_{\beta}
\subset \C^{N}$ \cor{(as shown in Figure \ref{analyticity:figure1}
  (a))}:
\begin{equation}
\Theta_{\beta} : = \left\{ \bz \in \mathbb{C}^{N};\,
\bz = \by + \bv,\,\by \in [-1,1]^{N},\,
 \sum_{n=1}^{N}  \sup_{x \in U}  \| \partial \bb_n \|_{2} 
\sqrt{\mu_{n}} |v_n| \cor{\leq \beta}
\right\}.
\label{analyticity:region}
\end{equation}

Now, we can extend the mapping $\partial F(\bmeta, \by) = I +
R(\bmeta,{\bf y})$, with \cor{$R(\bmeta,\by) := \sum_{n=1}^{N} $ $
  \sqrt{\mu_{n}} \partial \bb_n(\bmeta) y_{n}$}, to $\C^{N}$ by simply
replacing $\by$ with $\bz \in \Theta_{\beta}$. It is clear \cor{due to
  linearity} that the entries of the maps $F$ and $\partial F$ are
holomorphic in $\C^{N}$. Moreover, denote by $\Psi \equiv
F(\Theta_{\beta})$ \cor{the image} of $F:\Theta_{\beta} \rightarrow
\Psi$.

  Since $\by \in [-1,1]^{N}$ then the matrix inverse of $\partial
  F(\by)$ can be written as $\partial F^{-1}(\by) = (I + R(\by) )^{-1}
  = I + \sum_{k=1}^{\infty} \cor{(-R(\by))^{k}}$. Furthermore, since
  $\beta < \tilde{\delta}$ then the holomorphic expansion of $\partial
  F^{-1}(\by)$ can be written as the series
\[
\partial F^{-1}(\bz) = (I + R(\bz) )^{-1} = I +
\sum_{k=1}^{\infty} \cor{(-R(\bz))^{k}}
\]
and is \cor{pointwise} convergent $\forall \bz \in \Theta_{\beta}$. It
follows that each entry of $\partial F(\bz)^{-1}$ is analytic for all
$\bz \in \Theta_{\beta}$.

Up to this point we have assumed that only the geometry is stochastic
but have made no assumptions on further randomness in the forcing
function, the boundary conditions or the initial condition in Problems
\ref{setup:Prob1} and \ref{setup:Prob2}. These terms \corg{can also be}
extended with respect to other stochastic spaces.

\begin{assumption}~
\begin{enumerate}[(a)]
\item Suppose that the $N_{\bfq}$ valued random vector
  \cor{$\bfq:=[f_1, \dots, f_{N_{\bfq}}]^{T}$} takes values on
  $\Gamma_{\bfq}:= \tilde \Gamma_1 \times \dots \times \tilde
  \Gamma_{N_{\bfq}}$ with the probability density $
  \rho_{\bfq}(\bfq):\Gamma_{N_{\bfq}} \rightarrow [0,+\infty)$. The
    domains $\tilde \Gamma_1, \dots, \tilde \Gamma_{N_{\bfq}}$ can be
    assumed to be closed intervals in $\R$.  Now, assume that the
    random vector $\bfq$ is independent of $\by$ and write the forcing
    function $f:\Dw \times \Gamma_{\bfq} \rightarrow \R$ as
\[
\cor{f(\bx,\bfq,t) = 
\sum_{n=1}^{N_{\bfq}} c_{n}(t,f_n)
\xi_{n}(\bx),}
\]
where for $n = 1, \dots, N_{\bfq}$, \cor{$c_{n}(t,\bfq)\in
L^{\infty}_{\rho_{\bfq}}(\Gamma_{\bfq})$} $\forall t \in \R^{+}$, and
$\xi_{n}:\Dw \rightarrow \R$. Since $\xi_{n}$ is defined on $\Dw$ we
can remap $f:\Dw \times \Gamma_{\bfq} \rightarrow \R$ with pullback
onto the reference domain as
\[
(f \circ F)(\bmeta,\bfq,\by,t) = \sum_{n = 1}^{N_{\bfq}} \cor{c_{n}(t,f_n)}
(\xi_{n} \circ F)(\bmeta,\by).
\]
\cor{We shall now make analytic extension assumptions of the
  coefficients $c_{n}(t,\bfq)$ and $\xi_{n}$ for $n = 1,
  \dots,N_{\bfq}$.  The coefficients $c_{n}(\cdot,\bfq):\Gamma_{\bfq}
  \rightarrow \R$ are defined over the domain $\Gamma_{\bfq}$. Since
  the solution $\hat u$ from Problem \ref{setup:Prob2} is dependent on
  the coefficient $c_{n}(t,\bfq)$ certain analyticity assumptions have
  to be made.  In particular, suppose there exists an analytic
  extension of $c_{n}(\cdot,\bfq)$ onto the set $\msF \subset
  \C^{N_{\bfq}}$, where $\Gamma_{N_{\bfq}} \subset \C^{N_{\bfq}}$ (See
  Figure \ref{analyticity:figure1} for a graphical
  representation). The size of the region $\msF$ will directly depend
  on the coefficients $c_n(\cdot,\bfq)$ on a case by case
  basis. Furthermore, for $n = 1, \dots, N_{\bfq}$ the following
  assumptions are made:}
\begin{itemize}
\item $(\xi_{n} \circ F)(\bmeta,\by)$
can be analytically extended on $\Theta_{\beta}$, $\Real (\xi_{n}
\circ F)(\bz) \in L^{2}(U), \Imag (\xi_{n} \circ F)(\bz) \in L^{2}(U)$
$\forall \bz \in \Theta_{\beta}$.
\item $\Real \partial_{z_n} (\xi_{n} \circ F)(\bz), \Imag
\partial_{z_n} (\xi_{n} \circ F)(\bz) \in L^{2}(U)$ where
$\partial_{z_n}$ refers to the the Wirtinger derivative along the
$n^{th}$ dimension.
\end{itemize}

\item The initial condition $(u_0 \circ F)(\bmeta,\by)$ has an
  analytic extension on $\Theta_{\beta}$. Moreover, it is assumed that
  $\Real{ (u_0 \circ F)(\bmeta,\bz) }, \Imag{ (u_0 \circ
    F)(\bmeta,\bz) } \in L^{2}(U)$ for all $\bz \in \Theta_{\beta}$.
\end{enumerate}
\label{analyticity:assumption2}
\end{assumption}

\begin{assumption} \corg{We make the following assumptions on the
    Neumann boundary conditions:}
    It is also assumed than $(g_N \circ F)(\bmeta,\by)$
      can be analytically extended on $\Theta_{\beta}$, and that
      $\Real (g_N \circ F)(\bz) \in L^{2}(\partial U), \Imag (g_N
      \circ F)(\bz) \in L^{2}(\partial U)$ $\forall \bz \in
      \Theta_{\beta}$. Furthermore, assume that $det ( \bJ_{\tau}^{T}$
      $\partial F(\bmeta, \bz) ^{T}\partial F(\bmeta,\bz)
      \bJ_{\tau})^{\frac{1}{2}}$ is analytic for all $\bz$ in some
      region $\mcC \subset \C^{N}$ for all $\tau \in \mcT$.
\label{analyticity:assumption4}
\end{assumption}

\begin{remark}
  \corg{Since $\partial F(\bmeta,\bz)$ is analytic everywhere then
    $s(\bmeta,z):=det ( \bJ_{\tau}^{T} \partial F(\bmeta,\bz)^{T}
     \partial F(\bmeta,\bz)$ $\bJ_{\tau})$ is analytic in $\C^{N}$. Thus
    $s(\bmeta,\bz)^{\frac{1}{2}}$is analytic if $\Real s(\bmeta,\bz) >
    0 $. The region $\mcC \subset \C^{N}$ can be synthesized by
    placing the restriction on $\Real s(\bmeta,\bz) > 0$. This can be
    achieved by placing restrictions on $\partial F(\bmeta,\bz)$ for
    all $\bz \in \mcC$. This is, however, a little involved and is
    left for a future publication. Thus, to simplify the rest of the
    discussion in this paper we assume that there exists a constant
    $\hat \beta$ such that $\beta \leq \hat \beta < \tilde \delta$ and
    $\mcC = \Theta_{\beta} \subset \Theta_{\hat \beta}$.}
\end{remark}

\begin{figure}
  \begin{center}
    \begin{tabular}{c c}
\begin{tikzpicture}
    \begin{scope}

    \draw [line width=2] (-1, 0) -- (1, 0) node [below left]  {$\Gamma$}; 
    \draw [->] (-2.7,0) -- (2.7,0) node [above left]  {$\R^{N}$};
    \draw [->] (0,-1.5) -- (0,1.5) node [below right] {$i \R^{N}$};
    \draw [->,>=latex,dashed] (1,0) -- (1.70710678,0.70710678) node [above right] {$\bv$};

    \draw (1,-3pt) -- (1,3pt)   node [above] {$1$};
    \draw (-1,-3pt) -- (-1,3pt) node [above] {$-1$};
    \end{scope}

    \node[shape=semicircle,rotate=270,fill={rgb:red,143;green,250;blue,143},
      semitransparent,inner sep=12.7pt, anchor=south, outer sep=0pt]
    at (1,0) (char) {};
    \node[shape=semicircle,rotate=90,fill={rgb:red,143;green,250;blue,143},
      semitransparent,inner sep=12.7pt, anchor=south, outer sep=0pt]
    at (-1,0) (char) {}; \path
    [fill={rgb:red,143;green,250;blue,143},semitransparent]
    (-1.001,-1) rectangle (1.001,1.001);

    \node [below right,darkgray] at (1.3,1.7) {\Large $\Theta_{\beta}$};
\end{tikzpicture}
&
\begin{tikzpicture}
    \begin{scope}

    \filldraw[fill={rgb:red,143;green,250;blue,143}, semitransparent]
    (0,0) ellipse (2 and 1);

    \draw [line width=2] (-1, 0) -- (1, 0) node [below left]  {$\Gamma_{\bfq}$};


    \draw [->] (-2.7, 0) -- (2.7, 0) node [above left]  {$\R^{N_{\bfq}}$};
    \draw [->] (0,-1.5) -- (0,1.5) node [below right] {$i \R^{N_{\bfq}}$};
    \draw (1,-3pt) -- (1,3pt)   node [above] {$$};
    \draw (-1,-3pt) -- (-1,3pt) node [above] {$$};
    \end{scope}
    
    \node [below right] at (1.50,1.7) {\Large{$\cor{\msF}$}}; 
\end{tikzpicture} \\
\cor{(a)} & \cor{(b)} 
\end{tabular}
\end{center}
  \caption{\cor{Graphical representation of the sets $\Gamma$ and
      $\Gamma_{\bfq}$.  (a) $\Theta_{\beta} \subset \C^{N}$ is the
      extension of the set $\Gamma$ with respect to the parameter
      $\beta$. (b) Extension of $\Gamma_{\bfq}$ into the region
      $\cor{\msF} \subset \C^{N_{\bfq}}$}.}
\label{analyticity:figure1}
\end{figure}

\cor{To show that an analytic extension of the solution to Problem
  \ref{setup:Prob2} exists certain assumptions on the diffusion
  coefficient $a(\bx)$ are made. This assumption is left quite general
  and should be checked on a case by case basis.}

\begin{assumption} Suppose that the diffusion coefficient
    $a(\bx):{\mathcal G} \rightarrow \R$ is a deterministic map defined
    over the domain ${\mathcal G} := \cup_{\omega \in \Omega} \Dw$.
    \corg{Furthermore, assume there exists an analytic extension of
    $a(\bx)$ such that if $\bx \in \Psi$
    then}
\begin{enumerate}[i)]
\item $a_{max} c \geq \Real a(\bx) \geq a_{min} c$,
\item $|\Imag a(\bx)| < a_{min}$,
\end{enumerate}
where $c = 1 / tan(c_1)$ and $\pi / 8 > c_1 > 0$.
\label{analyticity:assumption3}
\end{assumption}

\corg{Let $G(\bz) := (a \circ F)(\bmeta,\bz) \partial F^{-1}(\bz)
  \partial F^{-T}(\bz)|\partial F(\bz)|$ for all $\bz \in
  \Theta_{\beta}$, we can now conclude that $G(\bz)$ is analytic for
  all $\bz \in \Theta_{\beta}$.}


The following lemma shows under what conditions the matrix $\Real
G(\bz)$ is positive definite and provides uniform bounds for the
minimum eigenvalue of $\Real G(\bz)$. This lemma is key to showing
that there exists an analytic extension of $\hat u(\bmeta,\by)$ on
$\Theta_{\beta}$.

\begin{lemma} Suppose
\[
0 < \beta < \min \{ \frac{\tilde{\delta} \log{\gamma_c}}{d +
  \log{\gamma_c}}, \sqrt{1 + \tilde{\delta}^2/2} - 1 \},
\]
where $\gamma_c := \frac{
  2 \tilde \delta^{d} +
  c(2-\tilde{\delta})^{d}}{\tilde{\delta}^{d} +
  c(2-\tilde{\delta})^{d}}$ then $\Real G(\bz)$ is positive definite
$\forall \bz \in \Theta_{\beta}$ and

\begin{enumerate}[(a)]
\item $\lambda_{min}(\Real G(\bz)^{-1}) \geq
  \mcA(\tilde{\delta},\beta,d,c_1,a_{min},a_{max}) 
  >0$ where
\[
\begin{split}
\mcA(\tilde{\delta},\beta,d,c_1,a_{max},a_{min}) 
&:=
\frac{  (2 - \tilde
  \delta)^{-d} (2 - \alpha(\beta))^{-1}
}
     {(a^2_{max}c^2 + a^2_{min})^{1/2}
     }
  \left(
cos\left(2c_1\right)
  \tilde{\delta}(\tilde{\delta} - 2 \beta)
  \right.
  \\
  &\left.
  -sin \left(2c_1\right)
  2 \beta(2 + (\beta - \tilde{\delta}))
  \right),
\end{split}
\]
and $\alpha(\beta):= 2 - \exp{\left(-\frac{d\beta}{\tilde \delta -
    \beta}\right)}$,
\item $\lambda_{max}(\Real G(\bz)^{-1}) \leq
  \mcR(\tilde{\delta},\beta,d,c_1,a_{min}) < \infty$ where
\[
\mcR(\tilde{\delta},\beta,d,c_1,a_{min}) 
:=
(a_{min}c)^{-1} \tilde \delta^{-d} \alpha(\beta)^{-1}
(
2 \beta ( 2 + \beta - \tilde \delta) 
+
(2 - \tilde{\delta} + \beta)^{2}).
\]
\item $\sigma_{max}(\Imag G(\bz)^{-1}) \leq \mcQ(\tilde{\delta},\beta,d,c_1,a_{min})
  < \infty$ where
  \[
\begin{split}
\mcQ(\tilde{\delta},\beta,d,c_1,a_{min}) 
&:=
(a_{min}c)^{-1} \tilde \delta^{-d} \alpha(\beta)^{-1}
(2 \beta(2 + (\beta - \tilde{\delta})) \\
&+
((2-\tilde{\delta}) + \beta)^{2} + \beta^{2}).
\end{split}
\]
\end{enumerate}
\label{analyticity:lemma3}
\end{lemma}
\begin{proof}
  \noindent
  {\bf (a)} From the proof in Lemma 5 in \cite{Castrillon2016} and
  Assumption \ref{analyticity:assumption1} we have that if $\beta <
  \tilde{\delta}/2$ then
\begin{equation}
\begin{split}
  \lambda_{min}(\Real \partial F(\bz)^{T} \partial F(\bz)) 
&\geq 
\tilde{\delta}(\tilde{\delta} - 2 \beta) 
> 0.
\end{split}
\label{analyticity:eqn5a}
\end{equation}
Furthermore, for all $\bz \in \Theta_{\beta}$,
\begin{equation}
\begin{split}
    \max_{i = 1,\dots,d}|\lambda_{i}(\Imag \partial F({\bf
  z})^{T} \partial F(\bz))|
&\leq  2 \beta(2 + (\beta - \tilde{\delta})),
\end{split}
\label{analyticity:eqn6}
\end{equation}
thus
\[
\begin{split}
\Real G(\bz)^{-1} 
&= \Real \Big(
\frac{(a_R(\bz) - i a_I(\bz))}{|a(\bz)|^2}
\frac{(\xi_R(\bz) - i \xi_I(\bz))}{|\xi(\bz)|^2} 
(\Real \partial F(\bz)^T \partial F(\bz) \\
&+
i \Imag \partial F(\bz)^T \partial F(\bz)
)
\Big)
\\
& =
\Real \Big(
\frac{e^{-i\theta_{a(\bz)}}}{|a(\bz)|} \frac{e^{-i\theta_{\xi(\bz)}}}{|\xi(\bz)|}
(\Real \partial F(\bz)^T \partial F(\bz)
+
i \Imag \partial F(\bz)^T \partial F(\bz)
)
\Big)
\end{split}
\]
\noindent where with a slight abuse of notation $\xi(\bz) :=
\xi_{R}(\bz) + i\xi_{I}(\bz) = |\xi(\bz)|e^{i\theta_{\xi(\bz)}} =
det(I + R({\bz}) )$ and $a(\bz):= |a(\bz)|e^{i\theta_{a(\bz)}} =
a_R(\bz) + ia_I(\bz) = \Real (a \circ F)(\bmeta,\bz) + i\Imag (a \circ
F)(\bmeta,\bz)$.

\cor{It is simple to check that $\Real \partial F(\bz)^T \partial
  F(\bz)$ and $\Imag \partial F(\bz)^T \partial F(\bz)$ are
  Hermitian. Let $\psi_{R}(\bz) := \Real a^{-1}(\bz)\xi^{-1}(\bz)$ and
  $\psi_{I}(\bz) := \Imag a^{-1}(\bz)\xi^{-1}(\bz)$. By applying the
  dual Lidskii inequality (if $A,B \in \mathbb{C}^{d \times d}$ are
  Hermitian then $\lambda_{min}(A+B) \geq \lambda_{min}(A) +
  \lambda_{min}(B)$) and assuming that $\psi_{R}(\bz) > 0$ it follows
  that} \cor{\begin{equation}
\begin{split}
\lambda_{min}(\Real G(\bz)^{-1}) 
&\geq 
\lambda_{min}(\psi_{R}(\bz)
(\Real
\partial F(\bz)^T \partial F(\bz)) 
- \psi_I(\bz) \Imag \partial F(\bz)^T \partial F(\bz)) \\
&\geq 
\lambda_{min}( \psi_{R}(\bz)
\Real
\partial F(\bz)^T \partial F(\bz)) +
\lambda_{min}( -\psi_I(\bz)
\Imag \partial F(\bz)^T \partial F(\bz)) \\
&\geq 
\psi_{R}(\bz) \lambda_{min}( 
\Real
\partial F(\bz)^T \partial F(\bz)) +
\lambda_{min}( -\psi_I(\bz)
\Imag \partial F(\bz)^T \partial F(\bz)) \\
&\geq 
\psi_{R}(\bz)\lambda_{min}(\Real
\partial F(\bz)^T \partial F(\bz)) \\
&- |\psi_I(\bz)|\max_{k =
  1,\dots,d}|\lambda_{k}(\Imag \partial F(\bz)^T \partial F(\bz))|.
\end{split}
\label{analyticity:eqn10}
\end{equation}}
\cor{The next step
  is place sufficient condition on $\xi(\bz)$, $a(\bz)$ and $\partial
  F(\bz)^{T} \partial F(\bz)$ such that equation
  \eqref{analyticity:eqn10} is greater than zero.}
\begin{enumerate}[I)]
\item \cor{First we determine for what range of values of $\beta$ the
  following inequality is satisfied:
  \begin{equation}
    \xi_R(\bz) \geq c|\xi_I(\bz)|
    \label{analyticity:eqn10a}
  \end{equation}
  for all $\bz \in \Theta_{\beta}$.  From Lemma 4 in
  \cite{Castrillon2016} $iii)$ we have that if $\alpha = 2 -
  \exp{\frac{d\beta}{\tilde \delta - \beta}} > 0$ then $\Real
  det(\partial F(\by))\geq \tilde \delta^{d} \alpha$ and $|\Imag
  det(\partial F(\by))| \leq (2 - \tilde \delta^{d})(1 - \alpha)$. Thus
  we need to solve for $\beta$ such that
  \[
  \xi_R(\bz) \geq
  \delta^{d} \alpha
  \geq
  c(2 - \tilde \delta^{d})(1 - \alpha)
  \geq c|\xi_I(\bz)|
  \]
  for all $\bz \in \Theta_{\beta}$. This is achieved if $\beta <
  \frac{\tilde{\delta} \log{ \gamma_c }}{d + \log{ \gamma_c }}$, where
  $\gamma_c := \frac{2 \tilde \delta ^{d} +
    c(2-\tilde{\delta})^{d}}{\tilde{\delta}^{d} +
    c(2-\tilde{\delta})^{d}}$.}

\item \cor{From Assumption \ref{analyticity:assumption3}} it follows
  that $a_R(\bz) > c|a_I(\bz)|$ if $\bz \in \Theta_{\beta}$.

\item \corg{From inequalities \eqref{analyticity:eqn5a}} and
  \eqref{analyticity:eqn6} it follows that if $\beta < \sqrt{1 +
    \tilde{\delta}^2/2} - 1$ then
\[
\lambda_{min}(\Real \partial F({\bf
    z})^{T} \partial F(\bz)) > \max_{k =
    1,\dots,d}|\lambda_{k}(\Imag \partial F(\bz)^{T} \partial
  F(\bz))|.
\]
\end{enumerate}
From I) - II) it follows that $\psi_{R}(\bz) > |\psi_{I}(\bz)|$ since
the angle of $\psi(\bz)$ is less than \cor{$\pi/4$} for all $\bz \in
\Theta_{\beta}$. However, an explicit expression can be derived:
\[
\psi_{R}(\bz) -
|\psi_{I}(\bz)| 
= |\psi(\bz)| (cos(\theta_{\psi(\bz)})
- sin(\theta_{\psi(\bz)})),
\]
where $|\psi(\bz)| = \frac{1}{|a(\bz)||\xi(\bz)|}$ and
\corj{$\theta_{\psi(\bz)} = -\theta_{a(\bz)} - \theta_{\xi(\bz)}$}.  
\corj{We observe from Assumption 7 that
\begin{align*}
\tan  \theta_{a(\bz)} = \frac{Im(a(\bz))}{Re(a(\bz))} < \frac{\lvert Im(a(\bz)) \rvert}{Re(a(\bz))} \\
\tan  (-\theta_{a(\bz)}) = \frac{-Im(a(\bz))}{Re(a(\bz))} < \frac{\lvert Im(a(\bz)) \rvert}{Re(a(\bz))}.
\end{align*}
\noindent It follows that $ \lvert \theta_{a(\bz)} \rvert < \frac{\pi}{8} $. Apply the same argument to $ \theta_{\xi(\bz)} $, we have $ \lvert \theta_{\xi(\bz)} \rvert < \frac{\pi}{8} $. It follow that
\begin{equation}
\theta_{\psi(\bz)} = - \theta_{a(\bz)} - \theta_{\xi(\bz)} \in (-\frac{\pi}{4}, \frac{\pi}{4}).
\end{equation}
\noindent Since $ \cos(\theta) > \sin(\theta), \forall \theta \in (-\frac{\pi}{4}, \frac{\pi}{4}) $, we obtain
\[
\psi_{R}(\bz) -
|\psi_{I}(\bz)| 
 > 0.
\]}
In particular, substituting equations \corg{\eqref{analyticity:eqn5a}
  and \eqref{analyticity:eqn6} in equation \eqref{analyticity:eqn10}}
we obtain that for all $\bz \in \Theta_{\beta}$
\[
\begin{split}
\lambda_{min}(\Real G(\bz)^{-1}) 
&\geq
\mcA(\tilde{\delta},\beta,d,c_1,a_{min},a_{max}) > 0.
\end{split}
\]
From London's Lemma \cite{london1981} it follows that $\Real G({\bf
  z})$ is positive definite $\forall \bz \in \Theta_{\beta}$.  

{\bf (b)} From the proof in Lemma 5 in \cite{Castrillon2016} and
Assumption \ref{analyticity:assumption1} we have that
\begin{equation}
\begin{split}
  \lambda_{max}(\Real \partial F(\bz)^{T} \partial F(\bz))  
&\leq
(2 - \tilde{\delta} + \beta)^{2}.
\end{split}
\label{analyticity:eqn12}
\end{equation}
From Assumption \ref{analyticity:assumption3} we have that
$|a(\bz)|^{-1} \leq (a_{min}c)^{-1}$ for all $\bz \in
\Theta_{\beta}$. From Lemma 4 in \cite{Castrillon2016}
$|\xi(\bz)|^{-1} \leq \tilde
\delta^{-d} \alpha(\beta)^{-1}$
for all $\bz \in \Theta_{\beta}$. We then have
that
\begin{equation}
  |\psi(\bz)| \leq (a_{min}c)^{-1} \tilde \delta^{-d} \alpha(\beta)^{-1}.
\label{analyticity:eqn13}
\end{equation}
Applying the Lidskii inequality (if $A,B \in \mathbb{C}^{d \times d}$
are Hermitian then $\lambda_{max}(A+B) \leq \lambda_{max}(A) +
\lambda_{max}(B)$) and \corg{substituting equations 
\eqref{analyticity:eqn5a}, \eqref{analyticity:eqn6},
\eqref{analyticity:eqn12} and \eqref{analyticity:eqn13}}
\[
\begin{split}
\lambda_{max}(\Real G(\bz)^{-1}) 
&\leq
|\psi_{R}(\bz)| 
\lambda_{max}(
\Real \partial F(\bz)^{T} \partial F(\bz) )  \\
&+
|\psi_{I}(\bz)| 
\max_{i}|\lambda_{i}( \Imag \partial F(\bz)^{T} \partial F(\bz) )|
 \\
&\leq
\frac{
\lambda_{max}(
\Real \partial F(\bz)^{T} \partial F(\bz) ) 
+ 
\max_{i}|\lambda_{i}( \Imag \partial F(\bz)^{T} \partial F(\bz) )|
}{|\psi(\bz) |^{-1}} \\
&\leq \mcR(\tilde{\delta},\beta,d,c_1,a_{min}) < \infty.
\end{split}
\]
  {\bf (c)} Similarly to (b), \cor{as shown in \cite{Castrillon2016},}
  it can be shown that
\begin{equation}
\begin{split}
\sigma_{max}(\Imag \partial F({\bf
  z})^{T} \partial F(\bz))
& \leq  2 \beta(2 + (\beta - \tilde{\delta})).
\end{split}
\label{analyticity:eqn14}
\end{equation}
and
\begin{equation}
\begin{split}
  \sigma_{max}(\Real \partial F(\bz)^{T} \partial F(\bz))   
&\leq
((2-\tilde{\delta}) + \beta)^{2} + \beta^{2}.
\end{split}
\label{analyticity:eqn15}
\end{equation}
From equations \eqref{analyticity:eqn13}, \eqref{analyticity:eqn14}
and \eqref{analyticity:eqn15} it follows that
\[
\begin{split}
\sigma_{max}(\Imag G(\bz)^{-1}) 
&\leq 
|\psi_{R}(\bz)|
\sigma_{max}(
\Imag \partial F(\bz)^{T} \partial F(\bz) )
\\
&+
|\psi_{I}(\bz)| 
\sigma_{max}( \Real \partial F(\bz)^{T} \partial F(\bz) )
\\
& \leq \mcQ(\tilde{\delta},\beta,d,c_1,a_{min}) < \infty.
\end{split}
\]
\end{proof}

\begin{lemma} 
\corb{For all  $\bz \in \Theta_{\beta}$ and $\bmeta \in U$
whenever}
\[
0 < \beta < min \{ \tilde{\delta} \frac{\log{\gamma_c}}{d +
  \log{\gamma_c}}, \sqrt{1 + \tilde{\delta}^2/2} - 1 \}
\]
then
\[
\lambda_{min}(\Real
  G(\bz)) \geq \varepsilon(\tilde{\delta},\beta,d,c_1,a_{max},a_{min}) > 0, 
\]
where $\varepsilon(\tilde{\delta},\beta,d,c_1,a_{max},a_{min})$ is
equal to
\[
\begin{split}
           \left(1 + \left(
             \frac{ \mcQ(\tilde{\delta},\beta,d,c_1,a_{min})}
                  {\mcA(\tilde{\delta},\beta,d,c_1,a_{min},a_{max})}
             \right) ^{2}
             \right)^{-1}\mcR(\tilde{\delta},\beta,d,c_1,a_{min})^{-1}.
\end{split}
           \]
\label{analyticity:lemma4}
\end{lemma}

\begin{proof} The proof essentially follows Lemma 6 in \cite{Castrillon2016}.
\end{proof}


The main result of this section can now be proven. 


\begin{theorem} Let $0 < \tilde{\delta} < 1$ then
  $\hat{u}(\bmeta,\by,\bfq,t)$ can be analytically extended on
  $\Theta_{\beta} \times \cor{\msF}$ if
\[
\beta < min \{ \tilde{\delta} \frac{\log{\gamma_c}}{d
  + \log{\gamma_c}}, \sqrt{1 + \tilde{\delta}^2/2} - 1 \}.
\]
\label{analyticity:theorem1}
\end{theorem}
\begin{proof}




\corg{Suppose that $\bV$ is a vector valued Hilbert space equipped with the
inner product $(\bgamma,\bv)_{\bV}$, where $\bv := [\vartheta_1\,
    \vartheta_2]^{T}$ and $\bgamma := [\gamma_1\, \gamma_2]^{T}$,
such that for all $\vartheta_1,\vartheta_2,\gamma_1,\gamma_2 \in V$
\[
{\bf (}\bgamma,\bv) := (\gamma_1,\vartheta_1) + (\nabla \gamma_1, \nabla
\vartheta_1) + (\gamma_2,\vartheta_2) + (\nabla \gamma_2, \nabla
\vartheta_2).
\]
Consider the extension of $(\by,\bfq) \rightarrow (\bz,\bq)$ on
$\Theta_{\beta} \times \cor{\msF}$. Let $\Phi(\by,\bfq,t) : =
\hat{u}(\by,\bfq,t)$ and consider the extension $\Phi = \Phi_R + i
\Phi_I$ on $\Theta_{\beta} \times \cor{\msF}$, where $\Phi_R : = \Real
\Phi$ and $\Phi_I := \Imag \Phi$. Let $\bPhi = [\Phi_{R}, \,
  \Phi_{I}]^{T}$, then the extension of $\Phi$ on $\Theta_{\beta}
\times \cor{\msF}$ is posed in the weak form as: Find $\bPhi \in
L^{2}(0,T;\bV)$, with $\partial_t \bPhi \in L^{2}(0,T;\bV^{*})$,
\corj{such that
\begin{equation}
  \begin{aligned}
 \int_U
  \partial_t \bPhi^{T}
  \bC(\bz)^{T}
\bv
+ \nabla \bPhi^T  \bG(\bz)^{T}
\nabla  \bv \,\mbox{d}\bmeta  
&=
\int_U  \hat \bfq(\bz,\bq,t) \cdot \bv \,\mbox{d}\bmeta  
&+
\sum_{\tau \in \mcT}
\int_{B^{0}_{r}}
\bg \cdot  \bv \,\mbox{d}\bx'
\\
 &
&\mbox{in $U \times (0,T)$} \\
\bPhi
&=\zzero &\mbox{on
    $\partial U_D \times (0,T)$} \\
\bPhi &= \bPhi_{0}
&\mbox{on
$U \times \{ t = 0 \}$}
  \end{aligned}
\label{analyticity:eqn2}
\end{equation}}
for all $\bv \in \bV$, \corj{where $\bv : = [\vartheta_1, \vartheta_2]^{T}$,}
\begin{align*}
 \bG(\bz)
 &:= \left(
\begin{array}{cc}
 G_{R}(\bz) &  -G_{I}(\bz) \\
 G_{I}(\bz) &   G_{R}(\bz) \\ 
\end{array} 
\right)
&
\hat \bfq(\bz,\bq,t)
&:= \left(
\begin{array}{c}
f_{R} \\
f_{I} \\
\end{array}
\right)
&
\bg(\bz)
&:= \left(
\begin{array}{c}
g^R_{N} \\
g^I_{N} \\
\end{array}
\right)
&
\0
&:= \left(
\begin{array}{c}
0\\
0\\
\end{array}
\right),
\\
\bC(\bz)
&:= \left(
\begin{array}{cc}
 c_{R}(\bz) &   -c_{I}(\bz) \\
 c_{I}(\bz) &    c_{R}(\bz) \\
\end{array} 
\right)
&
\bd(\bz)
&:= \left(
\begin{array}{c}
d_R \\
d_I \\
\end{array}
\right)
&
\bPhi_0(\bz)
&:= \left(
\begin{array}{c}
u^R_0\\
u^I_0\\
\end{array}
\right),
\end{align*}
$G_{R}(\bz):= \Real\{ G(\bz)\}$, $G_{I}(\bz):=\Imag \{ G(\bz) \}$,
$c_R(\bz) := \Real \{ |\partial F(\bz)| \}$, $c_I(\bz) := \Imag
\{|\partial F(\bz)| \}$, $f_R := \Real \{ (f \circ F)(\bq,\bz,t)
  |\partial F(\bz)|\}$, $f_I := \Imag \{ (f \circ F)(\bq,\bz,t)
  |\partial F(\bz)|\}$, $u^R_0 = \Real{(u \circ F)(\bz)}$, $u^I_0 =
\Imag{(u \circ F)(\bz)}$, $d_R(\bz) := \Real \{\nabla \cdot
  G(\bz) \nabla \hat \bw\}$, $d_I(\bz) := \Imag \{\nabla \cdot
  G(\bz) \hat \nabla \bw \}$,
  $g^R_N = \Real \{
(g_N \circ F)(\bmeta \circ \xi_{\tau}
,\bz)
det ( \bJ_{\tau}^{T} \partial F(\bmeta \circ \xi_{\tau},
  \bz)
  ^{T} \partial F(\bmeta \circ \xi_{\tau},\bz)
  \bJ_{\tau})^{\frac{1}{2}} \}
  $
  and
  $g^I_N = \Imag \{
(g_N \circ F)(\bmeta \circ \xi_{\tau}
,\bz)$
 $det ( \bJ_{\tau}^{T} \partial F(\bmeta \circ \xi_{\tau},
  \bz)
  ^{T} \partial F(\bmeta \circ \xi_{\tau},\bz)
  \bJ_{\tau})^{\frac{1}{2}} \}
  $
}
The system of equations (\ref{analyticity:eqn2}) has a unique solution
if $G_{R}$ is uniformly positive definite ($\lambda_{min}(G_{R}({\bf
  z}))$ $> 0$) since this implies that $\lambda_{min}( \bG(\bz)) > 0$
uniformly. From \cor{Lemma \ref{analyticity:lemma3}} this condition is
satisfied if $\bz \in \Theta_{\beta}$. Moreover, $\Phi(\bz,\bq,t)$
coincides with $\Phi(\by,\bfq,t)$ whenever $\bz \in \Gamma$ and $\bq
\in \Gamma_{\bfq}$ thus making it a valid extension of
$\Phi(\by,\bfq,t)$ on $\Theta_{\beta} \times \cor{\msF}$.

\cor{We now analyze the analytic regularity of the solution
  $\Phi(\bz,\bq,t)$ with respect to variables in $\bz$. However, it is
  not necessary to perform the analysis with respect to all the
  variables $\bz$ jointly. It is sufficient to show that
  $\Phi(\bz,\bq,t)$ is analytic with respect to each variable $z_n$,
  $n = 1,\dots,N$, separately. As shown at the end of the proof it can
  be concluded that $\Phi(\bz,\bq,t)$ is analytic in $\Theta_{\beta}
  \times \msF$.}

\cor{First, we concentrate on the $z_n$ variable of the vector
  $\bz$. Let $s= \Real z_n$ and $w = \Imag z_n$.  The first step is to
  show that the derivatives $\partial_{s} \Phi$ and $\partial_{w}
  \Phi$ exist on $\Theta_{\beta} \times \msF$.}  \corb{Consider the
  following weak problems:}

\begin{enumerate}[(a)]

\item \corj{Find $\phi \in L^{2}(0,T;\bV)$},
  with \corj{$\partial_t \phi \in L^{2}(0,T;\bV^{*})$, s.t.}
\begin{equation}
  \begin{split}
  & \corj{\int_U
  \partial_t \phi^{T} \bC(\bz)^{T} \bv
+ \nabla \partial_w  \bPhi^T  \bG(\bz)^{T}
\nabla  \bv \,\mbox{d}\bmeta } 
=
\int_U
  (-\partial_t
  \bPhi^{T}
  \partial_w \bC(\bz)^{T}
  \bv \,-
  \\
  &
   \nabla 
\bPhi^T  \partial_w \bG(\bz)^{T}
\nabla  \bv 
+
\partial_w \hat \bfq(\bz,\bq,t) \cdot \bv)
\, \mbox{d}\bmeta
\corg{+
\sum_{\tau \in \mcT}
\int_{B^{0}_{r}}
\partial_w
\bg \cdot  \bv \,\mbox{d}\bx'}
  \end{split}
\label{analyticity:eqn3}
  \end{equation}
  in $U \times (0,T)$ for all $\bv \in \bV$ and
\[
\begin{aligned}
 \corj{ \phi}
&=\zzero & &\mbox{(on
    $\partial U_D \times (0,T)$)} \\
\corj{\phi} &= \partial_w \bPhi_{0} 
& &\mbox{(on
$U \times \{ t = 0 \}$)} .
  \end{aligned}
\]

\item
  \corj{Find $\varphi \in L^{2}(0,T;\bV)$}, with
  \corj{$\partial_t \varphi \in L^{2}(0,T;\bV^{*})$, s.t.}
\begin{equation}
  \begin{split}
  & \int_U
  \corj{\partial_t
  \varphi^{T}
  \bC(\bz)^{T}
\bv
+ \nabla \varphi^T  \bG(\bz)^{T}
\nabla  \bv \,\mbox{d}\bmeta  }
=
\int_U
  (-\partial_t
  \bPhi^{T}
  \partial_s \bC(\bz)^{T}
  \bv \, -
  \\
  &
  \nabla 
\bPhi^T  \partial_s \bG(\bz)^{T}
\nabla  \bv 
+
\partial_s \hat \bfq(\bz,\bq,t) \cdot \bv) \,\mbox{d}\bmeta
+
\partial_s \hat \bd(\bz) \cdot \bv
\corg{+
\sum_{\tau \in \mcT}
\int_{B^{0}_{r}}
\partial_s
\bg \cdot  \bv \,\mbox{d}\bx'}
  \end{split}
\label{analyticity:eqn4}
  \end{equation}
  in $U \times (0,T)$ for all $\bv \in \bV$ and
\[
\begin{aligned}
  \corj{\varphi}
&=\zzero & &\mbox{(on
    $\partial U_D \times (0,T)$)} \\
\corj{\varphi} &= \partial_s \bPhi_{0}
& &\mbox{(on
$U \times \{ t = 0 \}$)}.
  \end{aligned}
\]
\end{enumerate}
Since $ \bG(\bz)$ is uniformly positive definite then
\eqref{analyticity:eqn3} - \eqref{analyticity:eqn4} have a unique
solution whenever $\bz \in \Theta_{\beta}$. \\

\corb{The next step is to integrate \eqref{analyticity:eqn3} with
  respect to $w$ and \eqref{analyticity:eqn4} with respect $s$.
  Now, since $\Theta_{\beta}$ is bounded and closed in $\C^{n}$ then
  the hypothesis of Fubini's Theorem holds. \corj{Furthermore, assume that
  $\phi $, $ \varphi$ and $ \bPhi $ belong in the space of smooth functions with compact support
   $\mathcal{C}_{c}^{\infty}(0,T;\bV^{*})$ with respect to the time variable. This
  is a reasonable assumption since $\mathcal{C}_{c}^{\infty}(0, T ; \bV^{*}) $ is
  dense in $L^{2}(0, T ; \bV^{*})$. By taking limits with respect to
  the norm of $L^{2}(0, T ; \bV^{*})$ in $\mathcal{C}_{c}^{\infty}(0,T;\bV^{*})$ }it
  follows that the order of integration and differentiation can be
  interchanged and thus}

\corj{
\begin{align*}
& \int_{U} \int (\partial_{t} \phi^{T} \bC(\bz)^{T} + \partial_{t} \bPhi^{T} \partial_{w} \bC(\bz)^{T}) \,\mbox{d} w \, \bv \,\mbox{d} \bmeta \\
+ & \int_{U} \int (\nabla \phi^{T} \bG(\bz)^{T} + \nabla \bPhi^{T} \partial_{w} \bG(\bz)^{T}) \,\mbox{d} w \nabla \bv \,\mbox{d} \bmeta =  \\
& \int_{U} \hat \bfq(\bz, \bq, t) \cdot \bv \,\mbox{d} \bmeta + \sum_{\tau \in \mathcal{T}} \int_{B_{\tau}^{0}} \bg \cdot \bv \,\mbox{d} \bx^{'} \\
& \int_{U} \int (\partial_{t} \varphi^{T} \bC(\bz)^{T} + \partial_{t} \bPhi^{T} \partial_{s} \bC(\bz)^{T}) \,\mbox{d} s \,\bv \,\mbox{d} \bmeta \\
+ & \int_{U} \int (\nabla \varphi^{T} \bG(\bz)^{T} + \nabla \bPhi^{T} \partial_{s} \bG(\bz)^{T}) \,\mbox{d} s \nabla \bv \,\mbox{d} \bmeta =  \\
& \int_{U} \hat \bfq(\bz, \bq, t) \cdot \bv \,\mbox{d} \bmeta + \sum_{\tau \in \mathcal{T}} \int_{B_{\tau}^{0}} \bg \cdot \bv \,\mbox{d} \bx^{'}. 
\end{align*}
}

\corb{ From equations \eqref{analyticity:eqn2},
  \eqref{analyticity:eqn3} and \eqref{analyticity:eqn4}, and since
  $\bv$ is arbitrary it follows that:}

\corj{
\begin{align}
\int (\partial_{t} \phi^{T} \bC(\bz)^{T} + \partial_{t} \bPhi^{T} \partial_{w} \bC(\bz)^{T}) \,\mbox{d} w & + \int (\nabla \phi^{T} \bG(\bz)^{T} + \nabla \bPhi^{T} \partial_{w} \bG(\bz)^{T}) \,\mbox{d} w \nonumber \\
= & \partial_{t} \bPhi^{T} \bC(\bz)^{T} + \nabla \bPhi^{T} \bG(\bz)^{T} ;
\label{analyticity:eqn3a}\\
\int (\partial_{t} \varphi^{T} \bC(\bz)^{T} + \partial_{t} \bPhi^{T}
\partial_{s} \bC(\bz)^{T}) \,\mbox{d} s & + \int (\nabla \varphi^{T}
\bG(\bz)^{T} + \nabla \bPhi^{T} \partial_{s} \bG(\bz)^{T}) \,\mbox{d}
s \nonumber \\ = & \partial_{t} \bPhi^{T} \bC(\bz)^{T} + \nabla
\bPhi^{T} \bG(\bz)^{T} \label{analyticity:eqn3b}.
\end{align}
}

\corj{We now claim that $ \bPhi = \int \phi \,\mbox{d} w$
  satisfies equation \eqref{analyticity:eqn3a} and
  $ \bPhi = \int
  \varphi \,\mbox{d} s $ satisfies equation \eqref{analyticity:eqn3b}
  . Following the same type of argument as above, we can interchange
  the order of differentiation and integration freely. By the
  fundamental theorem of calculus, we have that}

\corj{
\begin{align*}
& \partial_{t} \bPhi^{T} \bC(\bz)^{T} + \nabla \bPhi^{T} \bG(\bz)^{T} =  \partial_{t} \int \phi^{T} \,\mbox{d} w \bC(\bz)^{T} + \nabla \int \phi^{T} \,\mbox{d} w^{'} \bG(\bz)^{T} \\
= & \int \partial_{w} (\int \partial_{t} \phi^{T} \,\mbox{d} w^{'} \bC(\bz)^{T}) \,\mbox{d} w + \int \partial_{w} (\nabla \phi^{T} \,\mbox{d} w^{'} \bG(\bz)^{T}) \,\mbox{d} w \\
= & \int (\partial_{t} \phi^{T} \bC(\bz)^{T} + \partial_{t} (\int \phi^{T} \,\mbox{d} w^{'}) \partial_{w} \bC(\bz)^{T}) \,\mbox{d} w \\
+ & \int (\nabla \phi^{T} \bG(\bz)^{T} + \nabla (\int \phi^{T} \,\mbox{d} w^{'}) \partial_{w} \bG(\bz)^{T}) \,\mbox{d} w \\
= & \int (\partial_{t} \phi^{T} \bC(\bz)^{T} + \partial_{t} \bPhi^{T} \partial_{w} \bC(\bz)^{T}) \,\mbox{d} w + \int (\nabla \phi^{T} \bG(\bz)^{T} + \nabla \bPhi^{T} \partial_{w} \bG(\bz)^{T}) \,\mbox{d} w.
\end{align*}
}

\corb{From equation \eqref{analyticity:eqn3a} it follows that $\bPhi =
  \int \phi \, \mbox{d} w$. Following the same argument it can be
  shown that $\bPhi = \int \phi \, \mbox{d} s$.}

  \corj{It follows that there exists two functions $ \int \phi
    \,\mbox{d} w $ and $ \int \varphi \,\mbox{d} s $ that solve
    equation \eqref{analyticity:eqn2}, meanwhile $ \phi $ solves
    \eqref{analyticity:eqn3} and $ \varphi $ solves
    \eqref{analyticity:eqn4}. By uniqueness, we must have $ \bPhi =
    \int \phi \,\mbox{d} w = \int \varphi \,\mbox{d} s $. Hence we
    conclude that:} \corj{
\begin{equation*}
\partial_{w} \bPhi = \phi, \partial_{s} \bPhi = \varphi.
\end{equation*}
}
The second step is now to show that the Cauchy-Riemann conditions are
satisfied.  Consider the two functions $P(\bz) := \partial_{s}
\Phi_{R}(\bz) - \partial_{w} \Phi_{I}(\bz)$ and $Q(\bz) :=
\partial_{w} \Phi_{R}(\bz) + \partial_{s} \Phi_{I}(\bz)$, $\bP :=
        [P(\bz),\,Q(\bz)]^{T}$. \corg{First, let us write out explicitly
        equation \eqref{analyticity:eqn4} for the first term:
        \begin{equation}
  \partial_t \partial_s \bPhi^{T} \bC(\bz)^{T} \bv =        
  (\partial_t \partial_s \Phi_R c_R -
  \partial_t \partial_s \Phi_I c_I)\vartheta_1
  +
  (\partial_t \partial_s \Phi_R c_I -
  \partial_t \partial_s \Phi_I c_R)\vartheta_2.
  \label{analyticity:eqn4a}
\end{equation}
Second, for equation \eqref{analyticity:eqn3} exchange $\vartheta_1$
with $\vartheta_2$, and $\vartheta_2$ with $-\vartheta_1$ (Note, that
this is valid since equations \eqref{analyticity:eqn2} and
\eqref{analyticity:eqn3} are satisfied for all $\bv \in \bV$), then
the first term can written explicitly as
\begin{equation}
(\partial_t \partial_w \Phi_R c_R -
  \partial_t \partial_w \Phi_I c_I)\vartheta_2
  -(\partial_t \partial_w \Phi_R c_I -
  \partial_t \partial_w \Phi_I c_R)\vartheta_1.
  \label{analyticity:eqn5}
\end{equation}
      Adding Equations
      \eqref{analyticity:eqn4a} and \eqref{analyticity:eqn5} we
      obtain
\[
  \partial_t \bP^{T} \bC(\bz)^{T} \bv.
  \]
  }
Following for the rest of the terms we obtain the following weak
problem: Find $\bP \in L^{2}(0,T;\bV)$, with $\partial_t \bP \in
L^{2}(0,T;\bV^{*})$, s.t.
\corj{
\[
  \begin{split}
  & \int_U
  \partial_t
  \bP^{T}
  \bC(\bz)^{T}
\bv
+ \nabla \bP^{T}  \bG(\bz)^{T}
\nabla  \bv \,\mbox{d}\bmeta  \\
&=
\int_U
(-
  \partial_t
  \bPhi^{T}
 \left[
\begin{array}{cc}
  \partial_s c_{R}(\bz) -  \partial_w c_{I}(\bz)
  &  \partial_s c_{I}(\bz) +  \partial_w c_{R}(\bz)
  \\
  -(\partial_s c_{I}(\bz) + \partial_w c_{R}(\bz))
  & \partial_s c_{R}(\bz) -  \partial_w c_{I}(\bz)
  \\
\end{array} 
\right]
  \bv
  \\
  &+
    \nabla
  \bPhi^{T}
 \left[
\begin{array}{cc}
  \partial_s G_{R}(\bz) - \partial_w G_{I}(\bz)
  &  \partial_s G_{I}(\bz) +  \partial_w G_{R}(\bz)
  \\
  -(\partial_s G_{I}(\bz) + \partial_w G_{R}(\bz))
  &   \corj{\partial_s G_{R}(\bz) - \partial_w G_{I}(\bz)}
  \\
\end{array} 
\right]
\bv \\
&+
    [\partial_s f_R(\bz,\bq,t) - \partial_w f_I(\bz,\bq,t)\,\,\,
      \corj{\partial_s f_I(\bz,\bq,t) + \partial_w f_R(\bz,\bq,t)}
    ]^{\corg{T}}
    \\
    \,\mbox{d}\bmeta \\
    &\corj{+
\sum_{\tau \in \mcT}
\int_{B^{0}_{r}}
[\partial_s g^R_N(\bz) - \partial_w g^I_N(\bz)\,\,\,
      \partial_s g^I_N(\bz) + \corb{\partial_w g^R_N(\bz)}
]^{T}
\cdot  \bv \,\mbox{d}\bx'}
  \end{split}
  \]
  }
  in $U \times (0,T)$ for all $\bv \in \bV$ and
\[
\begin{aligned}
  \bP
&=\zzero & &\mbox{(on
    $\partial U_D \times (0,T)$ and $U \times \{ t = 0 \}$)}.
  \end{aligned}
\]
Since $(f \circ F)(\bq,\bz,t)$ is holomorphic in $\Theta_{\beta} \times
\cor{\msF}$ and $c(\bz)$ and $G(\bz)$ are holomorphic in
$\Theta_{\beta}$ then from the Cauchy Riemann equations we have that
\[
  \int_U
  \partial_t
  \bP^{T}
  \bC(\bz)^{T}
\bv
+ \nabla \bP^{T}  \bG(\bz)^{T}
\nabla  \bv \,\mbox{d}\bmeta  
= 0.
\]
\corb{Observe that zero solved the above equation
  above, and hence due to uniqueness} we have that $Q(\bz) = P(\bz) =
\zzero$ and therefore $\Phi(\bz,\bq,t)$ is holomorphic in
$\Theta_{\beta}$ along the $n^{th}$ dimension. From Hartog's Theorem
(Chap1, p32, \cite{Krantz1992}) and Osgood's Lemma (Chap 1, p 2,
\cite{Gunning1965}) $\Phi(\bz,\bq,t)$ is holomorphic in
$\Theta_{\beta}$ whenever $\bq \in \cor{\msF}$.

Since $\hat \bfq(\bz,\bq,t)$ is holomorphic in $\Theta_{\beta} \times
\cor{\msF}$ then $\Phi(\bz,\bq,t)$ is also holomorphic in
$\cor{\msF}$ whenever $\bz \in \Theta_{\beta}$. Applying Hartog's
Theorem and Osgood's Lemma it follows that $\Phi(\bz,\bq,t)$ is
holomorphic in \cor{$\Theta_{\beta} \times \cor{\msF}$}.

\end{proof}

\section{Stochastic polynomial approximation}
\label{stochasticcollocation}

Consider the problem of approximating a function $\nu: \Gamma
\rightarrow W$ on the domain $\Gamma$.  Our goal is to seek an
accurate approximation of $\nu$ in a suitably defined finite
dimensional space. To this end the following spaces are defined:

\begin{enumerate}[i)]

\item Let ${ \mcP_{\pp}(\Gamma)}
\subset L^2_\rho(\Gamma)$ be the span of tensor product polynomials of
degree at most $\pp = (p_1,\ldots,p_{N})$; i.e., $ \mcP_{\pp}(\Gamma)
= \bigotimes_{n=1}^{N}\;\mcP_{p_n}(\Gamma_{n})$ with $\mcP_{
  p_n}(\Gamma_{n}):=\text{\rm span}(y_n^m,\,m=0,\dots,p_n),$ $\quad
n=1,\dots,N$. 

\end{enumerate}


Suppose that $l^{\pp}_{k}$, $k \in {\mathcal K}$, is a series of Lagrange
polynomials that form a basis for ${ \mcP_{\pp}(\Gamma)}$.  An
approximation of $\nu$, know as the Tensor Product (TP)
representation, can be constructed as
\[
\nu^{N}(\by) = \sum_{k \in {\mathcal K}} \nu(\cdot, \cor{\mathbcal{y}_k})
l^{\pp}_{k}(\by)
\]
where \cor{$\mathbcal{y}_k$} are evaluation points from an appropriate
set of abscissas.  However, this is a poor choice for approximating
$\nu$ as the dimensionality of the index set ${\mathcal K}$ is $\Pi_{n =
  1}^{N} (p_n + 1)$. Thus the computational burden quickly becomes
prohibitive as the number of dimensions $N$ increases. This motivates
us to choose a reduced polynomial basis while retaining good accuracy.

Consider the univariate Lagrange interpolant along the $n^{th}$
dimension of $\Gamma$:
\[
{\mathcal I}^{m(i)}_{n}:C^{0}(\Gamma_n) \rightarrow {\mathcal P}_{m(i)-1}(\Gamma_n).
\]
In the above equation $i \geq 0$ is the level of approximation and
$m(i) \in \N_{0}$ is the number of evaluation points at level $i \in
\N_{0}$ where $m(0) = 0$, $m(1) = 1$ and $m(i) \leq m(i+1)$ if $i \geq
1$. \cor{Note that by convention $\mcP_{-1} =\emptyset$.}

An interpolant can now be constructed by taking tensor products of
${\mathcal I}^{m(i)}_n$ along each dimension $n$. However, the
dimensionality of $\mcP_p$ increases as $\prod_{n=1}^N$ $(p_n+1)$ with
$N$.  Thus even for a moderate size of dimensions the computational
cost of the Lagrange approximation becomes intractable. In contrast,
given \cor{sufficient regularity} of $\nu$ with respect to the random
variables defined on $\Gamma$, the application of Smolyak sparse grids
\cor{is better suited}
\cite{Smolyak63,Novak_Ritter_00,Back2011,nobile2008a}).

Consider the difference operator along the \cor{$n^{th}$ dimension of
  $\Gamma$}
\[
  {\Delta_n^{m(i)} :=} \mcI_n^{m(i)}-\mcI_n^{m(i-1)}.
\]
Given an integer $w \geq 0$, called the approximation level, and a
multi-index $\ii=(i_1,\ldots,i_{N})$ $\in \Nset^{N}_+$, let
$g:\Nset^{N}_+\rightarrow\Nset$ be a strictly increasing function in
each argument.

We can now construct a sparse grid from a tensor product of the
difference operators along every dimension.  However, the function $g$
imposes a restriction along each dimension such that a small subset of
the polynomial tensor is selected. More precisely, the sparse grid
approximation of $\nu$ is constructed as
\[
  \mcS^{m,g}_w[\nu]
= \sum_{\ii\in\Nset^{N}_+: g(\ii)\leq w} \;\;
 \bigotimes_{n=1}^{N} {\Delta_n^{m(i_n)}}(\nu(\by)) 
\]
\noindent or equivalently written as
\[
  \mcS^{m,g}_w[\nu(\by)]
  = \sum_{\ii\in\Nset^{N}_+: g(\ii)\leq w} \;c(\ii)\; 
  \bigotimes_{n=1}^{N} \mcI_n^{m(i_n)}(\nu(\by)), \,\,\, \text{with } c(\ii) 
  = \sum_{\stackrel{\jj \in \{0,1\}^{N}:}{g(\ii+\jj)\leq w}} (-1)^{|\jj|}.
\]

Let $\mm(\ii) = (m(i_1),\ldots,m(i_{N}))$ and consider the set of
polynomial multi-degrees
\[
\Lambda^{m,g}(w) = \{\pp\in\Nset^{N}, \;\;  g(\mm^{-1}(\pp+\oone))\leq w\}. 
\]
Denote by $\Pol_{\Lambda^{m,g}(w)}(\Gamma)$ the corresponding
multivariate polynomial space spanned by the monomials with
multi-degree in $\Lambda^{m,g}(w)$, i.e.
\[
\Pol_{\Lambda^{m,g}(w)}(\Gamma) = span\left\{\prod_{n=1}^{N} y_n^{p_n},
  \;\; \text{with } \pp\in\Lambda^{m,g}(w)\right\}.
\]


We have different choices for $m$ and $g$. One of the objectives is to
achieve good accuracy while restricting the growth of dimensionality
of the space $\Pol_{\Lambda^{m,g}(w)}(\Gamma)$.  A good choice of $m$
and $g$ is given by
\[
m(i) = \begin{cases} 1, & \text{for } i=1 \\ 2^{i-1}+1, & \text{for }
  i>1 \end{cases}\quad \text{ and } \quad g(\ii) = \sum_{n=1}^N
(i_n-1).
\]
For this choice the index set $
\cor{\Lambda^{m,g}(w)}:=\{\pp\in\Nset^{N}: \;\; \sum_n f(p_n) \leq
w\}$ where
\[
    f(p) = \begin{cases}
      0, \; p=0 \\
      1, \; p=1 \\
      \lceil \log_2(p) \rceil, \; p\geq 2
    \end{cases}.
\]

\begin{figure}[h]
\begin{center}
  \includegraphics[width=4in,height=2in, trim = {0 0 0 4.5cm},
    clip]{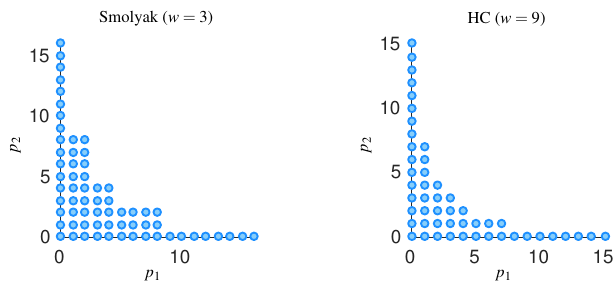}
\end{center}
\cor{\caption{Index sets for Smolyak (SM) sparse grid for $N = 2$ and $w =
  3$.  The Hyperbolic Cross (HC) index set is also shown for $N = 2$
  and $w = 9$, see \cite{Castrillon2016} for details.}}
\label{multivariate:fig1}
\end{figure}

This selection is known as the Smolyak sparse grid.  \cor{Other
  choices include the Total Degree (TD) and Hyperbolic Cross (HC),
  which are described in \cite{Castrillon2016}.  See Figure
  \ref{multivariate:fig1} for a graphical representation of the index
  sets $\Lambda_{m,g}(w)$ for $N = 2$.}
 

The Smolyak sparse grid combined with Clenshaw-Curtis (extrema of
Chebyshev polynomials) abscissas leads to nested sequences of one
dimensional interpolation formulas and a sparse grid with a highly
reduced number of nodes compared to the corresponding tensor grid. For
any choice of $m(i) > 1$ the Clenshaw-Curtis abscissas are given by
\[
\mathbcal{y}^n_j = -cos \left( \frac{\pi(j-1)}{m(i) - 1} \right).
\]
It is also straightforward to build related anisotropic sparse
approximation formulas by making the function $g$ to act differently
on the input random variables $y_n$ for $n = 1, \dots, N$. Anisotropic
sparse grids have been developed in \cor{\cite{Schillings2013,nobile2008b}}.

\section{Error analysis}
\label{erroranalysis}

In this section error estimates of the mean and variance of the QoI
are derived with respect to the sparse grid approximation and the
truncation of the stochastic model to the first $N_{s}$
dimensions. The error contributions from the finite element and
implicit solvers are neglected since there are many methods that can
be used to solve the parabolic equation (e.g. \cite{Knabner2003}) and
the analysis can be easily adapted. First, we establish some notation
and assumptions:

\begin{enumerate}[i)]

\item Split the Jacobian \cor{matrix} as follows
\begin{equation}
  \partial F(\bmeta,\omega) =
  I + \sum_{l=1}^{N_s} \sqrt{\mu_{l}}
  \partial \bb_{l}(\bmeta) Y_{l}(\omega)
  + \sum_{l=N_s+1}^{N}
  \sqrt{\mu_{l}} \partial \bb_{l}(\bmeta) Y_{l}(\omega).
\end{equation}
and
let $\Gamma_s := [-1,1]^{N_s}$, $ \cor{\Gamma_{\kappa}} :=
[-1,1]^{N-N_s}$, then the domain {$\Gamma = \Gamma_s \times
  \cor{\Gamma_{\kappa}}$}.
\cor{\item In practice one is interested in computing the statistics of a
  Quantity of Interest (QoI) of the solution over the stochastic
  domain or a subdomain of it. Assume that $Q:L^{2}(U) \rightarrow \R$
  is a bounded linear functional on $L^{2}(U)$ with norm $\| \cdot
  \|$.
\item Refer to $Q(\by_s)$ as $Q(\by)$ restricted to the stochastic
  domain $\Gamma_{s}$ and similarly for $G(\by_s)$. It is clear also
  that $Q(\by_s,\cor{\by_{\kappa}})=Q(\by)$ and
  $G(\by_s,\cor{\by_{\kappa}})=G(\by)$ for all {$\by \in \Gamma_s
    \times \cor{\Gamma_{\kappa}}$}, $\by_s \in \Gamma_s$, and
  $\cor{\by_{\kappa}} \in \cor{\Gamma_{\kappa}}$.
\item
Suppose that the $N_{\bg} < N_{\bfq}$ valued random vector
$\bg=[f_1,\dots,f_{N_{\bg}}]$ matches with $\bfq$ from the first to
$N_{\bfq}$ entry and takes values on $\Gamma_{\bg}:= \tilde \Gamma_1
\times \dots \times \tilde \Gamma_{N_{\bg}}$. The truncated forcing
function can now be written as
\[
(f \circ F)(\bmeta,\bg,\by,t) = \sum_{n = 1}^{N_{\bg}} c_{n}(t,f_n)
(\xi_{n} \circ F)(\bmeta,\by).
\]
}
\end{enumerate}

It is not hard to show that the \cor{variance error}
($|var[Q(\by_s,\cor{\by_{\kappa}},\bfq,t)] - var[\mcS^{m,g}_w[Q$
    $(\by_{s},\cor{\bg}, t)]] |$) and mean error
($|\eset{Q(\by_s,\cor{\by_{\kappa}},\bfq,t)} -
\eset{\mcS^{m,g}_w[Q(\by_s,\cor{\bg},t)}]|$) are less or equal to (see
\cite{Castrillon2016})
\[
\begin{split}
&C_{TR} \underbrace{\| Q(\by_{s},\cor{\by_{\kappa}},\bfq,t)
 - Q(\by_{s},\bfq,t)
 \|_{L^{2}_{\rho}(\Gamma \times \Gamma_{\bfq})}}_{\mbox{Truncation (I)}}
\\
&\cor{+ C_{FTR}
\underbrace{\| Q(\by_{s},\bfq,t)
 - Q(\by_{s},\bg,t)
 \|_{L^{2}_{\rho}(\Gamma \times \Gamma_{\bfq})}}_{\mbox{Forcing function Truncation (II)}}}
\\
&+
C_{SG}
\underbrace{
\| Q(\by_{s},\cor{\bg},t) - \mcS^{m,g}_w[Q(\by_{s},\cor{\bg},t)] 
\|_{L^{2}_{\rho}(\Gamma_s \times \Gamma_{\cor{\bg}})}
}
_{\mbox{\cor{Sparse Grid (III)}}},
\end{split}
\]
where $C_{TR}$, \cor{$C_{FTR}$} and $C_{SG}$ are positive constants
and $t \in (0,T)$. We now derive error estimates for the truncation
(I) and sparse grid (II) errors.

\subsection{Truncation error (I)}
\label{errorestimates:truncation}
Given that $Q: L^{2}(U) \rightarrow \mathbb{R}$ is a bounded linear
functional then
\[
\begin{split}
|Q(\by_s,\cor{\by_{\kappa}},\bfq,t) - Q(\by_s,\cor{\bfq},t)|
&\leq
\|Q\|
\| \hat{u}(\by_s,\cor{\by_{\kappa}},\bfq,t) - \hat{u}({\bf
  y}_s,\cor{\bfq},t) \|_{L^{2}(U)}.
\end{split}
\]
It follows that for $t \in (0,T)$
\[
\begin{split}
&
  \|Q(\by_s,\cor{\by_{\kappa}},\bfq,t) - Q(\by_s,\bfq,t)\|_{L^2_{\rho}(\Gamma \times
\Gamma_{\bfq}
  )} \\
&\leq \|Q\|
\| \hat{u}(\by_s,\cor{\by_{\kappa}},\bfq,t) - \hat{u}({\bf
  y}_s,\bfq,t) \|_{L^{2}_{\rho}(\Gamma \times \Gamma_{\bfq}; L^{2}(U))}.
\end{split}
\]
The objective now is to control the error term $\| \hat{u}(\by,\bfq,t)
- \hat{u}(\by_s,\bfq,t) \|_{L^{2}_{\rho}(\Gamma \times \Gamma_{\bfq};
  L^2(U))}$. But first we establish some notation.  If $W$ is a Banach
space defined on $U$ then let
\[
C^{0}(\Gamma; W) : = \{ v: \Gamma \rightarrow W\,\, \mbox{is
  continuous on $\Gamma$ and } \cor{\max_{\by \in \Gamma}
  \|v(\by)\|_{W}} < \infty \}.
\]
and
\[
L^{2}_{\rho}(\Gamma;W) := \{v:\Gamma \rightarrow W\,\, \mbox{is
  strongly measurable and} \,\, \int_{\Gamma}
\|v\|^{2}_{W}\,\rho(\by)\,\mbox{d}\by < \infty \}.
\]


With a slight abuse of notation let $\corg{\hat \varsigma}
(\by_s,\bfq,t) := \hat
u(\by_s,\bfq,t)$ for all $t \in (0,T)$, $\by_s \in \Gamma_{s}$ and
$\bfq \in \Gamma_{\bfq}$.  From Theorem \ref{analyticity:theorem1} it
follows that
\[
\corg{\hat \varsigma}, \hat{u} \in C^{0}(\Gamma \times \Gamma_{\bfq} ;L^{2}(0,T; V))
\subset L^2_{\rho}(\Gamma \times \Gamma_{\bfq};L^{2}(0,T; V)).
\]

\cor{The following theorem provides error bounds on the truncation
  error. It is adapted from Theorem 10 in \cite{Castrillon2016}.}
\begin{theorem} Suppose that
$\corg{\hat \varsigma} \in C^{0}(\Gamma_s;L^{2}(0,T; V))$
satisfies
\begin{equation}
\int_U
| \partial F(\by_s)
|
v \partial_t \corg{\hat \varsigma}\,
\mbox{\normalfont d}\bmeta
+ 
\cor{B(\by_s;\corg{\hat \varsigma},v)}
= 
\hat{l}(\by_s;\bfq,v)\,\,\,
\forall v \in V
\label{truncation:eqn1}
\end{equation}
for all $\bfq \in \Gamma_{\bfq}$, where $\corg{\hat \varsigma}(\by_s,\bfq,0) =
u_0$. Let $e(\by,\bfq,t):=\hat{u}(\by,\bfq,t) -
\corg{\hat \varsigma}(\by_s,\bfq,t)$,
  \[
  B_{\T} := \sup_{\bmeta \in U} \sum_{l = N_s+1}^{N} \sqrt{\mu_{l}}
  \| \partial \bb_{l}\|\,\,\mbox{and}\,\,C_{\T} := \sum_{i = N_s+1}^{N}
  \sqrt{\mu_{l}}\cor{\|\bb_l\|_{[L^{\infty}(U)]^{d}}  }
\]
  then for $0 < t < T$, $\bfq \in \Gamma_{\bfq}$, it follows that
\[
\begin{split}
\| e(\by,\bfq,t) \|^{2}_{L^2_{\rho}(\Gamma \times \Gamma_{\bfq}; L^2(U)) } 
&\leq \mathbb{C}_1 B_{\T}
+
\mathbb{C}_2 C_{\T}(1 + C_{\T}),
\end{split}
\]
\corj{where}
\corg{
\[
\begin{split}
&\mathbb{C}_1
(
T, C_{\mcT}, D_{\mcT},
C_{T}(U),C_P(U),
\F_{max}, \F_{min}, \tilde \delta, d, a_{max},
\| g_N \|_{W^{1,\infty}(\mcG_N)},
\\
& 
\sup_{t \in (0,T)} \eset{\| \hat u(\by,\bfq,t) \|^{2}_{V}},
\sup_{\small{{\begin{array}{c} t \in (0,T) \\ \bfq \in \Gamma_{\bfq},
\by \in \Gamma
\end{array}}}}
\| (f \circ F)(\by,\bfq,t) \|_{L^2(U)},
\\
&
\sup_{t \in (0,T)} \eset{
\| \partial_t \corg{\hat \varsigma}(\by_s,\bfq,t)\|^{2}_{L^{2}(U)}}^{\frac{1}{2}}
, \eset{\| \hat u(\by,\bfq,t) \|^{2}_{V}}^{\frac{1}{2}},
\sup_{\bx' \in B^{0}_{r}, \tau \in \mcT}\|\bJ_{\tau}(\bx')\|) \\
&\mathbb{C}_2
(T,S_{\mcT},C_{\mcT},C_T(U),C_P(U),\F_{max},\F_{min},d,
\sup_{\bfq \in \Gamma_{\bfq}} \|f\|_{W^{1,\infty}({\mathcal G} \times (0,T))}
, \|a\|_{W^{1,\infty}({\mathcal G})},
\\
&
\|g_N\|_{W^{1,\infty}(\partial {\mathcal G}_N)},
\| \chi_{U} \|_{L^{2}(U)}, 
\\
&
\sup_{t \in (0,T)} \corg{\eset{\| \hat u(\by,\bfq,t) \|^{2}_{V}}}
  ,\|
u_0
\|_{W^{1,\infty}({\mathcal G})})
\end{split}
\]
} \corg{are constants, $C_{T}(U)$ is the Trace Theorem constant,
  $C_{P}(U)$ is the Poincar\'{e} constants, $C_{\mcT} := (\inf_{\bx'
    \in B^{0}_{r}, \tau \in \mcT}$ $\sigma^{d-1}_{min}(\bJ_{\tau}^T
  \bJ_{\tau}))^{-1}$, $S_{\mcT} := \sup_{\bx' \in B^{0}_{r}, \tau \in
    \mcT, \by \in \Gamma} |s((\bmeta \circ \xi_{\tau})(\bx'),\by)^{\frac{1}{2}}|$ ,
  and $D_{\mcT} := (\inf_{\bx' \in B^{0}_{r}, \tau \in \mcT, \by \in
    \Gamma} |s((\bmeta \circ \xi_{\tau})(\bx'),\by)^{\frac{1}{2}}|
  )^{-1}$.}
\corg{\[
\chi_{U}(\bmeta) = \left\{
\begin{array}{l l c}
  1 & & \bmeta \in U \\
  0 & & o.w.
\end{array}
\right.
.
\]}
\label{errorestimates:theorem1}
\end{theorem}
\begin{proof}
Consider the solution to equation \eqref{truncation:eqn1}
\[
 \corg{\hat \varsigma} \in C^{0}(\Gamma_s \times \Gamma_{\bfq};L^{2}(0,T; V))
 \subset L^2_{\rho}(\Gamma \times \Gamma_{\bfq}; L^{2}(0,T; V))
\]
where the matrix of coefficients $G(\by_s)$ depends only on the
variables $Y_{1}, \dots, Y_{N_{s}}$. \corg{By adapting the proof from
Strang's Lemma we have that}
\[
\begin{split}
  \| \corg{\hat \varsigma}(\by_s) 
- \hat{u}(\by) \|^{2}_{V}
&\leq
\Kset \Big(
\left|
\hat{l}(\by_s;\corg{\hat \varsigma}(\by_s) - \hat{u}(\by)) - 
\hat{l}(\by;\corg{\hat \varsigma}(\by_s) - \hat{u}(\by))
\right|
\\
&+ 
\int_U 
(\corg{\hat \varsigma}(\by_s) - \hat{u}(\by))
(
| \partial F(\by) |
-
| \partial F(\by_s) |
)
\partial_t \corg{\hat \varsigma}(\by_s)
)
\\
&+ 
\int_U 
(\corg{\hat \varsigma}(\by_s) - \hat{u}(\by))
(
| \partial F(\by) |
(
\partial_t \hat u(\by)
-
\partial_t \corg{\hat \varsigma}(\by_s)
)
)
\\
&+ 
\left|
B(\by;\hat{u}(\by),\corg{\hat \varsigma}(\by_s) - \hat{u}(\by)) 
-
B(\by_s;\hat{u}(\by),\corg{\hat \varsigma}(\by_s) - \hat{u}(\by)) 
)
\right|
\Big)
,
\end{split}
\]
\corg{where $\Kset := a_{min}^{-1}\F_{min}^{-d} \F_{max}^{2}(1 +
C_{P}(U)^2)$.}
Recall that
$e(\by) : = \hat{u}(\by) - \corg{\hat \varsigma}(\by_s)$ and note that
\[
\begin{split}
\int_U e(\by)|\partial F(\by)|^{\frac{1}{2}}
\partial_t \cor{
  \left( |\partial F(\by)|^{\frac{1}{2}} e(\by) \right)
}
&=
\frac{1}{2}
\partial_t
\|e(\by) |\partial F(\by)|^{\frac{1}{2}}\|^{2}_{L^2(U)} \\
&\geq
\frac{\F_{min}^{d}}{2}
\partial_t \| e(\by) \|^{2}_{L^2(U)}
\end{split}
\]
thus
\[
\begin{split} 
\frac{\F^{d}_{min}}{2} 
\partial_t \| e(\by) \|^{2}_{L^2(U)}
&\leq
\Big(
\left|
\hat{l}(\by;e(\by)) - 
\hat{l}(\by_s;e(\by))
\right| \\
&+ \left|
B(\by;\hat{u}(\by),e(\by)) 
-
B(\by_s;\hat{u}(\by),e(\by)) 
)
\right|
\\
&+ 
\int_U 
\left|
e(\by)
(
| \partial F(\by) |
-
| \partial F(\by_s) |
)
\partial_t \corg{\hat \varsigma}(\by_s)
\right|
\Big).
\end{split}
\]
and for all $t \in (0,T)$, $\bfq \in \Gamma_{\bfq}$
and $\by \in \Gamma$
\[
\begin{split}
\partial_t \| e(\by, \bfq,t) \|^{2}_{L^2(U)}
&\leq 
\frac{2}{\F^{d}_{min}}
(\Bset_1
+
\Bset_2
+
\Bset_3
)
\end{split}
\]
for some non-negative constants $\Bset_1,\Bset_2,\Bset_3 <
\infty$. For now assume that $\Bset_1$, $\Bset_2$ and $\Bset_3$ are
known. From Gronwall's inequality we have that for $t \in (0,T)$, $\by
\in \Gamma$, and $\bfq \in \cor{\Gamma_{\bfq}}$
\begin{equation}
\begin{split}
\| e(\by,\bfq,t) 
\|^{2}_{L^{2}(U)} 
&\leq 
\| e(\by,\bfq,0) 
\|^{2}_{L^{2}(U)} 
+ 
\frac{2(\Bset_1
+
\Bset_2
+
\Bset_3
)
T}{\F_{min}^{d}}
\end{split}
\label{erroranalysis:eqn4}
\end{equation}
The first term in equation \eqref{erroranalysis:eqn4} is bounded as
\begin{equation}
\begin{split}
\| e(\by,\bfq,0) \|_{L^{2}(U)} 
&= 
\| (u_0 \circ F)(\by_s)
-
(u_0 \circ F)(\by)
\|_{L^{2}(U)} \\
&\leq
\|
u_0
\|_{W^{1,\infty}({\mathcal G})}
\| \chi_U \|_{L^{2}(U)}
\sup_{\by \in \Gamma, \bmeta \in U} |F(\by_s) - F(\by)|,
\end{split}
\label{erroranalysis:eqn9}
\end{equation}
for all $\bfq \in \Gamma_{\bfq}$ and $\by \in \Gamma$. For the second
term we have that
\corg{
\[
\begin{split} 
&\Bset_1
:= \cor{\sup_{t \in (0,T)}}
| B(\by; \hat u(\by,\bfq,t), e(\by,\bfq,t)) -
   B(\by_s; \hat u(\by,\bfq,t), e(\by,\bfq,t))| \\
 & \leq 
\cor{\sup_{t \in (0,T)}}  
\| \hat{u}(\by,\bfq,t) \|_{V}
(
\| \hat{u}(\by,\bfq,t) \|_{V}
+
\| \corg{\hat \varsigma}(\by_s,\bfq,t) \|_{V}
)
\sup_{\bmeta \in U, \by \in \Gamma}  
\|G(\by) - G(\by_s)\|.
\end{split}
\]
}
Following the same argument for Theorem 10 in \cite{Castrillon2016} we
have that
\begin{equation}
\sup_{\bmeta \in U, \by \in \Gamma}  
\|G(\by) - G(\by_s)\| \leq a_{max}B_{\T}H(\F_{max},\F_{min},\tilde \delta,d)
\label{erroranalysis:eqn5}
\end{equation}
for some constant $H(\F_{max},\F_{min},\tilde \delta,d)$. Thus we have 
\begin{equation}
  \begin{split}
    \eset{\Bset_1}
    &\leq
    a_{max}B_{\T}H(\F_{max},\F_{min},\tilde \delta,d) 
    \cor{\sup_{t \in (0,T)}
    2
    \eset{\| \hat u(\by,\bfq,t) \|^{2}}.}
\end{split}
\label{erroranalysis:eqn6}
\end{equation}
\cor{The last term is true since $\hat{u}$ and $\corg{\hat \varsigma}$ are equal
  when $\by = [\by_s, 0]^{T}$
  \[
\eset{\| \corg{\hat \varsigma}(\by_s,\bfq,t) \|^{p}_{V}}
\leq
\eset{\| \hat{u}(\by,\bfq,t) \|^{p}_{V}},
\]
$p = 1, 2, \dots$ where the expectation $\eset{\cdot}$ is defined over the
domain $\Gamma$ and $\Gamma_{\bfq}$.} The next constant
\[
\Bset_2 :=
|\hat{l}(\by;e(\by,\bfq,t) ) - 
\hat{l}(\by_s;e(\by,\bfq,t) )|
\]
\corj{is bounded by}
\begin{equation}
\begin{split}
& 
\Big|
\int_U ( (f \circ F)(\by,\bfq,t) |\partial F(\by)|  
- 
(f \circ F)(\by_s,\bfq,t) |\partial F(\by_s)|)
e(\by,\bfq,t)
\Big|
\\
&+
\corg{\Big| \sum_{\tau \in \mcT}
  \int_{B^{0}_{r}}
  (
(g_N \circ F)(\bmeta
,\by)
s(\bmeta,\by)^{\frac{1}{2}}
  -
  (g_N \circ F)(\bmeta
,\by_s)
s(\bmeta,\by_s)^{\frac{1}{2}}
  )e(\by,\bfq,t)\,
\mbox{d}\bx' \Big|}
\\
& \leq 
\corg{
  \sum_{\tau \in \mcT}
  \int_{B^{0}_{r}}
  |(
(g_N \circ F)(\bmeta
,\by)
  -
  (g_N \circ F)(\bmeta
,\by_s))
s(\bmeta,\by))^{\frac{1}{2}}
e(\by,\bfq,t)|
}
\\
&
\corg{+
| 
(g_N \circ F)(\bmeta
,\by_s)
(s(\bmeta,\by)^{\frac{1}{2}}
  -
s(\bmeta,\by_s)^{\frac{1}{2}})
  e(\by,\bfq,t)|\,
  \mbox{d}\bx'}
\\
&
+\int_U 
|((f \circ F)(\by,\bfq,t)  - (f \circ F)(\by_s,\bfq,t))
|\partial F(\by)|e(\by,\bfq,t)|
\\
&+ 
|(f \circ F)(\by_s,\bfq,t)(|\partial F(\by)| 
- |\partial F(\by_s)|)e(\by,\bfq,t)| 
\\
\end{split}
\label{erroranalysis:eqn3}
\end{equation}
for all $t \in (0,T)$, $\bfq \in \Gamma_{\bfq}$ and $\by \in \Gamma$.
\corj{The following inequalities are used to bound equation}
\eqref{erroranalysis:eqn3}:
\corg{
  \begin{subequations}
  \label{erroranalysis:inequalities}
  \begin{eqnarray}
    &
    \sum_{\tau \in \mcT}
  \int_{B^{0}_{r}}
  |(
(g_N \circ F)(\bmeta
,\by)
  -
  (g_N \circ F)(\bmeta
,\by_s))
s(\bmeta,\by)^{\frac{1}{2}})
e(\by,\bfq,t)|
     \\
     &\leq
     C_{\mcT} S_{\mcT} \| e(\by,\bfq,t) \|_{L^{2}(\partial U)}
\nonumber 
\|g_N\|_{W^{1,\infty}(\partial {\mathcal G}_N)}
\sup_{\by \in \Gamma, \bmeta \in U}|F(\by) - F(\by_s)| \\
&(\mbox{Using the Trace Theorem \cite{Evans1998}}:
\| e(\by,\bfq,t) \|_{L^{2}(\partial U)}
  \leq
  C_{T}(U)\| e(\by,\bfq,t) \|_{H^{1}(U)}
  ),
\nonumber
\end{eqnarray}
\begin{eqnarray}
  &
   \sum_{\tau \in \mcT}
  \int_{B^{0}_{r}}  
| 
(g_N \circ F)(\bmeta
,\by_s)
(s(\bmeta,\by)^{\frac{1}{2}}
  -
s(\bmeta,\by_s)^{\frac{1}{2}})
  e(\by,\bfq,t)|
\\
&\leq
\frac{D_{\mcT}C_{\mcT}
C_{T}(U)}{2} \| e(\by,\bfq,t) \|_{H^{1}(U)}
\|g_N\|_{W^{1,\infty}(\partial {\mathcal G}_N)} 
\sup_{\tau \in \mcT}
|det(\bJ_{\tau}^{T}
\partial F(\bmeta \circ \xi_{\tau},\by)^{T}
\nonumber \\
&\partial F(\bmeta \circ \xi_{\tau},\by)
\bJ_{\tau}) 
-
det(\bJ_{\tau}^{T}
\partial F(\bmeta \circ \xi_{\tau},\by_s)^{T}
\partial F(\bmeta \circ \xi_{\tau},\by_s)
\bJ_{\tau})
|
\nonumber 
\end{eqnarray}
\begin{eqnarray}
&\int_U 
|((f \circ F)(\by,\bfq,t)  - (f \circ F)(\by_s,\bfq,t))
|\partial F(\by)| e(\by,\bfq,t) \\
&\leq
\F_{max}^{d} \| \chi_{U} \|_{L^{2}(U)}
\sup_{\bfq \in \Gamma_{\bfq}}
\|f\|_{W^{1,\infty}({\mathcal G} \times (0,T))}
\nonumber \\
&
\| e(\by,\bfq,t) \|_{V}
\sup_{\by \in \Gamma, \bmeta \in U}|F(\by) - F(\by_s)|, \nonumber
\end{eqnarray}
\begin{eqnarray}
& 
\int_U
|(f \circ F)(\by_s,\bfq,t)(|\partial F(\by)| 
- |\partial F(\by_s)|)e(\by,\bfq,t)| \\
&\leq
\|e(\by,\bfq,t)\|_{V}
\sup_{\small{{\begin{array}{c} t \in (0,T) \\ \bfq \in \Gamma_{\bfq},
\by \in \Gamma
\end{array}}}}
\| (f \circ F)(\by,\bfq,t) \|_{L^2(U)}
\nonumber   \\
&
\sup_{\by \in \Gamma, \bmeta \in U} 
||\partial F(\by)| 
- |\partial F(\by_s)||, \nonumber 
\end{eqnarray}
\label{erroranalysis:inequalities2}
  \end{subequations}
  }

\noindent
\corg{Following the same argument for Theorem 10 in \cite{Castrillon2016} we
have that
\begin{eqnarray}
  \label{erroranalysis:eqn7}
\sup_{\by \in \Gamma} 
||\partial F(\by)| 
- |\partial F(\by_s)||
\leq \F_{max}^{d-1}\F^{-2}_{min}dB_{\T}, \\
\sup_{\bmeta \in U, \by \in \Gamma} 
|F(\by)
- F(\by_s)| \leq C_{\T}.
\end{eqnarray}
From Theorem 2.12 in \cite{Ipsen2008} ($A,E \in \C^{d \times d}$ then
$|det(A + E) - det(A)| \leq d \|E\|$ $\max \{ \|A\|,$ $\|A +
E\|\}^{d-1})$ we obtain $\forall x\in U$ and $\forall \by \in \Gamma$
\begin{equation}
\begin{split}
&|det(\bJ_{\tau}^{T}
\partial F(\bmeta \circ \xi_{\tau},\by)^{T}
\partial F(\bmeta \circ \xi_{\tau},\by)
\bJ_{\tau}) 
-
det(\bJ_{\tau}^{T}
\partial F(\bmeta \circ \xi_{\tau},\by_s)^{T} \\
&\partial F(\bmeta \circ \xi_{\tau},\by_s)
\bJ_{\tau})|
\leq
\sup_{\bx' \in B^{0}_{r}, \tau \in \mcT}\|\bJ_{\tau}(\bx')\|^{2d}
d\F_{max}^{2d-1}B_{\T}.
\end{split}
\end{equation}
}
Furthermore using Jensen's inequality
\begin{equation}
\eset{\|e(\by,\bfq,t) \|_{V}} \leq C
\eset{\|\hat{u}(\by,\bfq,t) \|_{V}}
\leq C \eset{\|\hat{u}(\by,\bfq,t) \|^{2}_{V}}
\label{erroranalysis:eqn8}
\end{equation}
for some constant $C > 0$.  \corj{Combining inequalities
\eqref{erroranalysis:inequalities} (a) - (g) and equations
\eqref{erroranalysis:eqn3},
\eqref{erroranalysis:eqn7}
- \eqref{erroranalysis:eqn8}}
\corg{
\[
\begin{split}
&\eset{\Bset_2}
\leq 
\Xi
(C_T(U),S_{\mcT},C_{\mcT},\F_{max},\F_{min},d,
\sup_{\bfq \in \Gamma_{\bfq}} \|f\|_{W^{1,\infty}({\mathcal G} \times (0,T))}
, \|a\|_{W^{1,\infty}({\mathcal G})}, \\
&
\|g_N\|_{W^{1,\infty}(\partial {\mathcal G}_N)},
\|\hat \bw\|_{W^{1,\infty}(\mcG)},
\| \chi_{U} \|_{L^{2}(U)}, \\
&
\sup_{t \in (0,T)}
\eset{\| \hat u(\by,\bfq,t) \|^2_{V}}
) C_{\T} +
\Upsilon(
C_{T}(U),C_{\mcT},D_{\mcT},
\F_{max}, \F_{min}, \tilde \delta, d, a_{max}, \\
&
\| \hat g_N \|_{W^{1,\infty}(\partial \mcG_N)},
\sup_{\bx' \in B^{0}_{r}, \tau \in \mcT}\|\bJ_{\tau}(\bx')\|,
\\
&
\sup_{t \in (0,T)}
\eset{\| \hat u(\by,\bfq,t) \|^2_{V}},
\sup_{\small{{\begin{array}{c} t \in (0,T) \\ \bfq \in \Gamma_{\bfq},
\by \in \Gamma
\end{array}}}}
\| (f \circ F)(\by,\bfq,t) \|_{L^2(U)}
)
B_{\T}, 
\end{split}
\]
}
for some non-negative constants $\Xi$ and $\Upsilon$. The last
constant
\[
\begin{split}
\Bset_{3}
&\leq
\int_U 
\left|
e(\by,\bfq,t)
(
| \partial F(\by) |
-
| \partial F(\by_s) |
)
\partial_t \corg{\hat \varsigma}(\by_s,\bfq)
\right| \\
&\leq
2
\F^{d-1}_{\max}\F^{-2}_{min}dB_{\T} \sup_{t \in (0,T)}
\| \hat u(\by,\bfq,t) \|_{V}
\| \partial_t \corg{\hat \varsigma}(\by_s,\bfq,t)\|_{L^{2}(U)}
.
\end{split}
\]
\corg{By using the Schwartz inequality $\eset{\Bset_3}$
is less or equal to
\[
\F^{d-1}_{\max}\F^{-2}_{min}dB_{\T} \sup_{t \in (0,T)}
\Big(\eset{\| \hat u(\by,\bfq,t)\|^{2}_{V}}\Big)^{1/2}\Big(\eset{
  \| \partial_t \corg{\hat \varsigma}(\cor{\by_s},\bfq,t)\|^{2}_{L^{2}(U)}}
\Big)^{1/2}.
\]
}
Combining the bounds for $\eset{\Bset_1}$ , $\eset{\Bset_2}$,
$\eset{\Bset_3}$, equations \eqref{erroranalysis:eqn9} and
\eqref{erroranalysis:eqn4} we obtain the result.
\end{proof}

\cor{
\subsection{Forcing function truncation error (II)} 
Since $Q$ is a bounded linear functional the error due to (II) is
controlled by $\| \hat{u}(\by_s,\bfq,t) - \hat{u}(\by_s,\bg,t)
\|_{L^{2}_{\rho}(\Gamma \times \Gamma_{\bfq}; L^2(U))}$.  Suppose that
$\hat{u}(\by_s, \bfq,t) \in L^{2}(0,T; V)$ satisfies the following
equation
  \begin{equation}
\int_U
| \partial F(\by_s)
|
v \partial_t \hat{u}\,
\mbox{\normalfont d}\bmeta
+ 
B(\by_s;\hat{u},v)
= 
\hat{l}(\by_s;\bfq,v)\,\,\,
\forall v \in V
\label{forcingtruncation:eqn1}
\end{equation}
for all $\bfq \in \Gamma_{\bfq}$ and $\by_s \in \Gamma_s$, where
$\hat{u}(\by_s,\bfq,0) = u_0$. Furthermore, let $\hat{u}(\by_s,
\bg,t) \in L^{2}(0,T; V))$ satisfies
  \begin{equation}
\int_U
| \partial F(\by_s)
|
v \partial_t \hat{u}\,
\mbox{\normalfont d}\bmeta
+ 
B(\by_s;\hat{u},v)
= 
\hat{l}(\by_s;\bg,v)\,\,\,
\forall v \in V
\label{forcingtruncation:eqn2}
\end{equation}
for all $\bg \in \Gamma_{\bg}$ and $\by_s \in \Gamma_s$, where
$\hat{u}(\by_s,\bg,0) = u_0$.
\begin{theorem}
Let $\hat{e}(\by_s,\bfq,t) : = \hat{u}(\by_s, \bfq,t) - \hat{u}(\by_s,
\bg,t)$, $t \in (0,T)$,
\[
0 < \epsilon < a_{min}^{-1} \F_{min}^{-d} \F^{2}_{max}C_{P}(U)^{2}/4
\]
and
\[
\mcJ(d,a_{min}, \F_{min},\F_{max}, C_P(U),\epsilon) := 
\frac{2}{\F_{min}^{d}}\left[
  \frac{1}{4\varepsilon} -
  a_{min}\F^{d}_{min} \F^{-2}_{max} C_P(U)^{-2}
  \right]
\]
then
\[
\begin{split}
&\|\hat{e}(\by_s,\bfq,t) \|_{L^2_{\rho}(\Gamma \times \Gamma_{\bfq};U)}
  \leq
T^{1/2}e^{\mcJ(d,a_{min}, \F_{min},\F_{max}, C_P(U),\epsilon)T/2}
\\
&\epsilon^{1/2}
\left(\sum_{n=N_{\bg}+1}^{N_{\bfq} }
\eset{c^{2}_{n}(t,f_n)}\right)^{1/2}
\left( \sum_{n=N_{\bg}+1}^{N_{\bfq} }
\|(\xi_{n} \circ F)(\bmeta,\by_s)
\|^{2}_{L^{2}_{\rho}(\Gamma_s;U)}
\right)^{1/2}
.
\end{split}
\]
\label{erroranalysis:theorem2}
\end{theorem}
\begin{proof}
Subtract
\eqref{forcingtruncation:eqn2} from \eqref{forcingtruncation:eqn1}
\begin{equation}
\int_U
| \partial F(\by_s)
|
v \partial_t \hat{e}\,
\mbox{\normalfont d}\bmeta
+ 
B(\by_s;\hat{e},v)
= 
\int_U ((f \circ F)(\cdot,\by_s,\bfq)
-
(f \circ F)(\cdot,\by_s,\bg))v
\label{forcingtruncation:eqn3}
  \end{equation}
  $\forall v \in V$. Recall that
\[
\begin{split}
\int_U \hat{e}|\partial F(\by_s)|^{\frac{1}{2}}
\partial_t 
  \left( |\partial F(\by_s)|^{\frac{1}{2}} \hat{e}\right)
&=
\frac{1}{2}
\partial_t
\|\hat{e} |\partial F(\by_s)|^{\frac{1}{2}}\|^{2}_{L^2(U)}. 
\end{split}
\]
Let $v = \hat{e}$ and substitute in \eqref{forcingtruncation:eqn3},
then
\[
\begin{split}
&  \frac{1}{2}
 \partial_t
 \|\hat{e} |\partial F(\by_s)|^{\frac{1}{2}}\|^{2}_{L^2(U)}
+ 
B(\by_s;\hat{e},\hat{e}) \\
&= 
\int_U ((f \circ F)(\cdot,\by_s,\bfq)
-
(f \circ F)(\cdot,\by_s,\bg))\hat{e}.
\end{split}
\]
Applying the Poincar\'{e} and Cauchy's inequalities we obtain
\[
\begin{split}
&  \frac{\F^{d}_{min}}{2}
 \partial_t
 \|\hat{e} \|^{2}_{L^2(U)}
+ a_{min}\F^{d}_{min} \F^{-2}_{max} C_P(U)^{-2}\|\hat{e}\|^{2} \\
&\leq
\frac{1}{4\epsilon}\| \hat{e} \|^2_{L^{2}(U)}
+
\epsilon
\|(f \circ F)(\cdot,\by_s,\bfq) -
(f \circ F)(\cdot,\by_s,\bg))\|^2_{L^{2}(U)}.
\end{split}
\]
From Gronwall's inequality it follows that
\[
\begin{split}
  \eset{\|\hat{e} \|^{2}_{L^2(U)}}
  &\leq
Te^{\mcJ(d,a_{min}, \F_{min},\F_{max}, C_P(U),\epsilon)T}
\\
&\epsilon
\eset{\|(f \circ F)(\cdot,\by_s,\bfq) -
(f \circ F)(\cdot,\by_s,\bg))\|^2_{L^{2}(U)}}.
\end{split}
\]
We have that
\[
\begin{split}
&\|(f \circ F)(\cdot,\by_s,\bfq) -
  (f \circ F)(\cdot,\by_s,\bg))\|_{L^{2}(U)} \\
&  \leq 
\| \sum_{n=N_{\bg}+1}^{N_{\bfq} } c_{n}(t,\bfq) (\xi_{n} \circ F)(\bmeta,\by_s)
\|_{L^{2}(U)} \\
&\leq 
\sum_{n=N_{\bg}+1}^{N_{\bfq} } |c_{n}(t,f_n)| \|(\xi_{n} \circ F)(\bmeta,\by_s)
\|_{L^{2}(U)} \\
&\leq
\Big(
\sum_{n=N_{\bg}+1}^{N_{\bfq} }
c^{2}_{n}(t,f_n)
\Big)^{1/2}
\Big(
\sum_{n=N_{\bg}+1}^{N_{\bfq} }
\|(\xi_{n} \circ F)(\bmeta,\by_s)
\|^{2}_{L^{2}(U)}
\Big)^{1/2},
\end{split}
\]
thus
\[
\begin{split}
& \eset{
\|(f \circ F)(\cdot,\by_s,\bfq) -
(f \circ F)(\cdot,\by_s,\bg)) \|^{2}_{L^{2}(U)}
} \\
&\leq
\sum_{n=N_{\bg}+1}^{N_{\bfq} }
\eset{c^{2}_{n}(t,f_n)}
\sum_{n=N_{\bg}+1}^{N_{\bfq} }
\|(\xi_{n} \circ F)(\bmeta,\by_s)
\|^{2}_{L^{2}_{\rho}(\Gamma_s;U)}.
\end{split}
\]
\end{proof}
}
\subsection{\cor{Sparse grid error (III)}} 
In this section convergence rates for the isotropic Smolyak sparse
grid with Clenshaw Curtis abscissas are derived. The \cor{convergence}
rates can be extended to a larger class of abscissas and anisotropic
sparse grids following the same approach.

Given the bounded linear functional $Q:L^{2}(U) \rightarrow \R$ it
follows that
\[
\begin{split}
|Q(\by_s,\cor{\bg},t) - \mcS^{m,g}_w[Q(\by_s,\cor{\bg},t)]|
&\leq
\|Q\|
\| \hat{u}(\by_s,\cor{\bg},t) - \mcS^{m,g}_w
   [\hat{u}(\by_s,\cor{\bg},t)] \|_{L^{2}(U)}
\end{split}
\]
for all $t \in (0,T)$, $\by_s \in \Gamma_s$ and $\cor{\bg} \in
\Gamma_{\cor{\bg}}$. The sparse grid operator $\mcS^{m,g}_w$ is with
respect to the domain $\Gamma_{s} \times \Gamma_{\cor{\bg}}$.  The next step
it to bound the term
\[
\| \hat{u}(\by_s,\cor{\bg},t) - \mcS^{m,g}_w [\hat{u}(\by_s,\cor{\bg},t)]
\|_{L^{2}(\Gamma_s \times \Gamma_{\cor{\bg}}; U)}.
\]
for $t \in (0,T)$.  The error term $\| \epsilon \|_{L^{2}(\Gamma_s
  \times \Gamma_{\cor{\bg}}; U)}$, where
\[
\epsilon :=
\hat{u}(\by_s,\cor{\bg},T)- \mcS^{m,g}_w [\hat{u}(\by_s,\cor{\bg},T)],
\]
is controlled by the number of collocation knots $\eta$ (or work), the
choice of the approximation formulas $(m(i), g({\bf i}))$, and the
region of analyticity of $\Theta_{\beta} \times \cor{\msF} \subset
\mathbb{C}^{N_s + N_{\cor{\bg}}}$. From Theorem \ref{analyticity:theorem1}
the solution $\hat{u}(\by_s,\cor{\bg},t)$ admits an analytic extension in
$\Theta_{\beta} \times \cor{\msF} \subset \C^{N_s+N_{\cor{\bg}}}$ for
all $t \in (0,T)$.


In \cite{nobile2008b,nobile2008a} the authors derive error estimates
for isotropic and anisotropic Smolyak sparse grids with
Clenshaw-Curtis and Gaussian abscissas where $\| \epsilon
\|_{{L^{2}_{\hat{\rho}}(\Gamma_{s};V})}$ exhibit algebraic or
sub-exponential convergence with respect to the number of collocation
knots $\eta$.
For these estimates to be valid the solution $\hat{u}(\by_s,\cor{\bg},T)$
has to admit and extension on a polyellipse in $\C^{N_s + N_{\cor{\bg}}}$
i.e. ${\mathcal E}_{\sigma_1, \dots, \sigma_{N_s + N_{\cor{\bg}}}} : =
\Pi_{i=1}^{N_s + N_{\cor{\bg}}}$ ${\mathcal E}_{n,\sigma_n}$, where
\[
\begin{split}
{\mathcal E}_{n,\sigma_n} 
&=
 \left\{ z \in \C;\,\Real(z) = \frac{e^{\sigma_n}
 +
  e^{-\sigma_n}}{2}cos(\theta), \right. \\
&
\left.
\Imag(z)= \frac{e^{\sigma_n} -
  e^{-\sigma_n}}{2}sin(\theta),\theta \in [0,2\pi) \right\},\,\,\,
\end{split}
\]
and $\sigma_n > 0$. For an isotropic sparse grid the overall
asymptotic subexponential decay rate $\hat{\sigma}$ will be dominated
by the smallest $\sigma_n$ i.e.
\[
\hat{\sigma} \equiv \min_{n = 1,\dots, N_s + N_{\cor{\bg}}} \sigma_n.
\]
Then the goal is to choose the largest $\hat{\sigma}$ such that ${\mathcal
  E}_{\sigma_1, \dots, \sigma_{N_s + N_{\cor{\bg}}}} \subset \Theta_{\beta}
\times \cor{\msF}$. First, form the set $\Sigma \subset \C^{N_s}$
such that $\Sigma \subset \Theta_{\beta}$, where $\Sigma:=\Sigma_1
\times \dots \times \Sigma_{N_s}$ and
\[
\Sigma_{n} := \left\{ \bz \in \mathbb{C};\, \bz = \by
+ \bv,\,\by \in [-1,1],\, |v_n| \leq \tau_n :=
\frac{\beta}{1 - \tilde{\delta}}
\right\},
\]
for $n = 1,\dots,N_s$. Let
\[
\hat \sigma_{\beta} := \log{ \left(\sqrt{ \left( \frac{\beta}{1 -
      \tilde \delta} \right)^{2} + 1} + \frac{\beta}{1 - \tilde
    \delta} \right)} > 0,
\]
then the polyellipse ${\mathcal E}_{\sigma_1, \dots, \sigma_{N_{\cor{\bg}}}}$ can be
embedded in $\Sigma$ by setting $\sigma_1 = \sigma_2 = \dots =
\sigma_{N_s} = \hat{\sigma}_{\beta}$, as shown in Figure
\ref{errorestimates:figure1}.

The second step is to form a polyellipse such that $ {\mathcal
  E}_{\sigma_1, \dots, \sigma_{N_{\cor{\bg}}}} \subset \cor{\msF}$ . This,
of course, depends on the size of the region $\cor{\msF}$. For
simplicity \cor{we assume that} $\sigma_{N_s+1} = \sigma_{N_s+2} =
\dots = \sigma_{N_s + N_{\cor{\bg}}} = \hat{\sigma}_{\cor{\bg}}$, for some
constant $\hat{\sigma}_{\cor{\bg}} > 0$. The constant $\hat{\sigma}_{\cor{\bg}}$
is chosen such that ${\mathcal E}_{\sigma_{N_s + 1}, \dots, \sigma_{N_s +
    N_{\cor{\bg}}}} \subset \cor{\msF}$. Finally, the polyellipse ${\mathcal
  E}_{\sigma_1, \dots, \sigma_{N_s + N_{\cor{\bg}}}}$ is embedded in
$\Theta_{\beta} \times \cor{\msF}$ by setting $\hat \sigma =
\min\{\sigma_{\beta}, \sigma_{\cor{\bg}}\}$.

From Theorem 3.11 \cite{nobile2008a}, given that $w > \frac{N_s +
  N_{\cor{\bg}}}{\log{2}}$ for a nested CC sparse grid the following estimate
holds:
\begin{equation}
  \begin{split}
    {\| \epsilon \|_{L^{2}_{\hat{\rho}}(\Gamma_{s} \times \Gamma_{\cor{\bg}};V)}}
&\leq
    {\mathcal
  Q}(\sigma,\delta^{*},N_s + N_{\cor{\bg}})\eta^{\mu_3(\sigma,\delta^{*},N_s + N_{\cor{\bg}})} \\
&* \exp
  \left(-\frac{(N_s + N_{\cor{\bg}}) \sigma}{2^{1/(N_s + N_{\cor{\bg}})}} \eta^{\mu_2(N_s + N_{\cor{\bg}})} \right)
\end{split}
  \label{erroranalysis:sparsegrid:estimate}
\end{equation}
where 
\[{\mathcal Q}(\sigma,\delta^{*}, \tilde N) := 
\frac{C_1(\sigma,\delta^{*})}{\exp(\sigma \delta^{*} \tilde{C}_2(\sigma)  )}
\frac{\max\{1,C_1(\sigma,\delta^{*})\}^{\tilde N}}{|1 - C_1(\sigma,\delta^{*})|}, 
\]
$\sigma = \hat{\sigma}/2$, $\tilde N \in \N_{+}$, $\mu_2(\tilde N) =
\frac{log(2)}{\tilde N(1 + log(2(\tilde N)))}$ and
$\mu_3(\sigma,\delta^{*},\tilde N) = \frac{\sigma
  \delta^{*}\tilde{C}_2(\sigma)}{1 + \log{(2(\tilde N))}}$.  The
  constants $C_1(\sigma,\delta^{*})$, $\tilde{C}_2(\sigma)$ and
  $\delta^{*}$ are defined in \cite{nobile2008a} equations (3.11) and
  (3.12).

\begin{figure}
\begin{center}
\begin{tikzpicture}
  \begin{scope}[thick,font=\scriptsize]


    \draw [->] (-2.5,0) -- (2.5,0) node [above left]  {$\Real $};
    \draw [->] (0,-2) -- (0,2) node [below right] {$\Imag$};
    \draw (1,-3pt) -- (1,3pt)   node [above] {$1$};
    \draw (-1,-3pt) -- (-1,3pt) node [above] {$-1$};
    \end{scope}

    \node[shape=semicircle,rotate=270,fill=lightgray,semitransparent,inner
      sep=12.7pt, anchor=south, outer sep=0pt] at (1,0) (char) {};
    \node[shape=semicircle,rotate=90,fill=lightgray,semitransparent,inner
      sep=12.7pt, anchor=south, outer sep=0pt] at (-1,0) (char) {};
    \path [draw=none,fill=lightgray,semitransparent] (-1.001,-1) rectangle
    (1.001,1.001);
    \filldraw[fill={rgb:red,143;green,250;blue,143},semitransparent]
    (0,0) ellipse (2 and 1);
    \node [below right,black] at
    (1.50,1.25) {$\Sigma_n \subset \Theta_{\beta}$}; \node [below
      right,darkgray] at (0.85,0.85) {${\mathcal E}_{\sigma_n}$};
\end{tikzpicture}
\end{center}
\caption{Embedding of the polyellipse ${\mathcal
    E}_{\sigma_1,\dots,\sigma_{N_s}}:=\Pi^{N_s}_{n = 1} {\mathcal
    E}_{n,\sigma_n}$ in $\Sigma \subset \Theta_{\beta}$. Each ellipse
  ${\mathcal E}_{n,\sigma_n}$ is embedded in $\Sigma_n \subset
  \Theta_{\beta}$ for $n = 1,\dots, N_s$ .}
\label{errorestimates:figure1}
\end{figure}

\begin{remark}
  \cor{Note that for the convergence rate given by equation
    \eqref{erroranalysis:sparsegrid:estimate} there is an implicit
    assumption that the constant $M(u(\bz_s, \bq,t)):= \max_{\bz_s \in
      \Theta_{\beta}, \bq \in \cor{\msF}} \|\hat u(\bz_s,
    \bq,t)\|_{V}$, for $t \in (0,T)$, is equal to one.  This
    assumption was introduced in \cite{nobile2008a} to simplify the
    overall presentation of the convergence results. This constant for
    $t \in (0,T)$ can be easily reintroduced in equation
    \eqref{erroranalysis:sparsegrid:estimate}.}  However, it will not
  change the overall convergence rate.
\end{remark}

\section{Numerical results}
\label{numericalresults}

\begin{figure}[h]
\hspace{-13mm}
\begin{tabular}{c}
\includegraphics[ width=5.2in, height=2.1in]{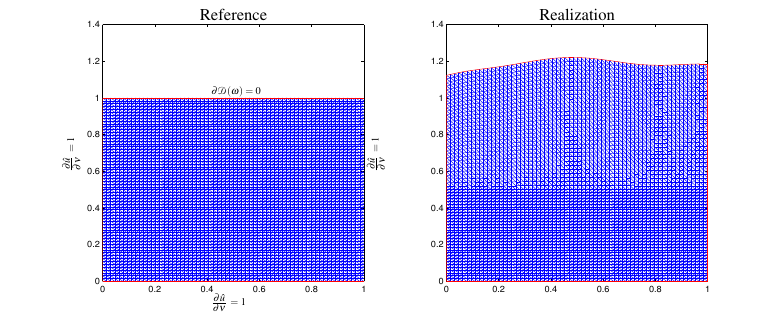}
\\
\hspace{8mm} (a) \hspace{53mm} (b) \\
\includegraphics[ width=2.5in, height=2.1in]{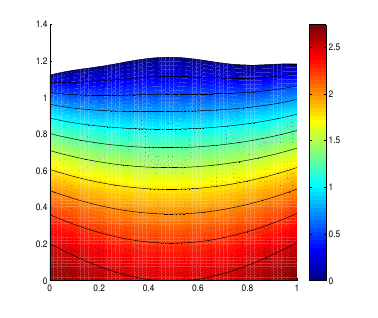}
\\
(c) \\
\end{tabular}
  \caption{Stochastic deformation of a square domain and solution on a
    realization of the stochastic domain.  (a) Reference square domain
    with Dirichlet boundary conditions. (b) Vertical deformation from
    stochastic model. (c) Contours of the solution of the parabolic
    PDE for $T = 1$ on the stochastic deformed domain realization.}
\label{numericalresults:fig1}
\end{figure}

In this section numerical examples are executed that elucidate
\cor{the} truncation and Smolyak sparse grid convergence rates for
parabolic PDEs.  Suppose the reference domain \cor{is set to $U :=
  (0,1) \times (0,1)$} and is deformed according to the following
rule:
\[
\begin{array}{llll}
  F(\eta_{1}, \eta_{2}) = (\eta_{1},\,(\eta_{2}-0.5)(1 +
  ce(\omega,\eta_{1})) + 0.5) & & if & \eta_{2} > 0.5\\ F(\eta_{1},
  \eta_{2}) = (\eta_{1},\,\eta_{2}) & & if & 0 \leq \eta_{2} \leq 0.5
\end{array}
\]
\noindent for some positive constant $c > 0$. This deformation rule
only stretches (or compresses) the upper half of the domain and fixes
the button half. The Dirichlet boundary conditions are set to zero for
the upper border. The rest of the borders are set to \cor{Neumann}
boundary conditions with $\frac{\partial u}{\partial \nu} = 1$ (See
Figure \ref{numericalresults:fig1} (a)). \cor{Furthermore, the
  diffusion coefficient $a(\bx) = 1$ and the forcing function $f =
  0$.}

The stochastic model $e(\omega, \eta_1)$ is defined as
\begin{center}
$e_{S}(\omega,\eta_{1}) := Y_{1}(\omega)\left(
    \frac{\sqrt{\pi}L}{2} \right)+ \sum_{n = 2}^{N_{s}} {\sqrt{\lambda_{n}}}
  \varphi_{n}(\eta_{1})Y_{n}(\omega); \hspace{1mm}$ $e_{F}(\omega,\eta_{1})
  := \sum_{n = N_{s}+1}^{N} {\sqrt{\lambda_{n}}}
  \varphi_{n}(\eta_{1})Y_{n}(\omega)$,
\end{center}
where $\{Y_{n}\}_{n = 1}^{N}$ are independent uniform distributed in
$(-\sqrt{3},\sqrt{3})$. \cor{Note that through a rescaling of the
  random variables $Y_1(\omega), \dots, Y_N(\omega)$ the random vector
  $\bY(\omega):=[Y_1(\omega),\dots,Y_N(\omega)]$ can take values on
  $\Gamma$. Thus the analyticity theorems and convergence rates
  derived in this article are valid.}

To make a comparison between the theoretical decay rates and the
numerical results the gradient terms $\sqrt{\lambda_{n}} \sup_{x \in
  U} \|B_n(x)\|$ are set to decay linearly as \cor{$n^{-k}$, where $k
  = 1$ or $k = 1/2$,} thus for $n = 1, \dots, N$ let
${\sqrt{\lambda_{n}}} := \frac{(\sqrt{\pi}L)^{1/2}}{n}$, $n \in \N$,
and
\[
\varphi_{n}(\eta_{1}) : = \left\{
\begin{array}{cc}
  n^{-1}sin \left( \frac{\lfloor n/2 \rfloor \pi \eta_{1}}{L_{p}} \right) &
  \mbox{if n is even}\\ n^{-1}cos \left( \frac{\lfloor n/2 \rfloor \pi
      \eta_{1}}{L_{p}} \right) & \mbox{if n is odd}\\
\end{array}
\right.
\]
With this choice ${\sup_{x \in U} \sigma_{max}(B_{n}(x))}$, for $n =
1, \dots, N$, is bounded by a constant, \cor{which depends on $N$},
and linear decay on the gradient of the deformation is obtained.

The QoI is defined on the bottom half of the reference domain, which
is not deformed, as
\[
Q(\hat{u}(\omega, T)) := \int_{(0,1)} \int_{(0,1/2)}
\varphi(\eta_{1})\varphi(2\eta_{2})\hat{u}(\eta_1,\eta_2,\omega,T)
\,\mbox{d}\eta_{1}\mbox{d}\eta_{2},
\]
where $\varphi(x):= \exp \left( \frac{-1}{1 - 4(x-0.5)^2} \right)$.
\cor{The chosen QoI $Q$ can, for example, represent the weighed
  total chemical concentration in the region defined by $(0,1) \times
  (0,1/2)$ given uncertainty in the region. Other useful applications
  include sub-surface aquifers with soil variability, heat transfer,
  etc.}

To solve the parabolic PDE a finite element semi-discrete
approximation is used for the spatial domain. For the time evolution
an implicit second order trapezoidal method with a step size of $t_d$
and final time $T$.

\corg{For each realization of the domain $\Dw$ the mesh is perturbed
  by the function deformation $F$.} In Figure
\ref{numericalresults:fig1} the original reference domain (a) is
shown. An example realization of the deformed domain from the
stochastic model and the contours of the solution for the final time
$T = 1$ are shown in Figure \ref{numericalresults:fig1} (a) \&
(b). Notice the significant deformation of the stochastic domain.



\begin{remark}
  \cor{For $N = 15$ dimensions, $k = 1$ and $k = 1/2$ the mean
    $\mathbb{E}[Q(\hat{u}(\by))]$ and variance $\var[Q(\hat{u}(\by))]$
    are computed with a dimensional adaptive sparse grid method
    collocation with \corg{$\approx 10,000$} collocation points and a
    Chebyshev abscissa \cite{Gerstner2003}. For the linear decay, $k =
    1$, the computed normalized mean value is 0.9846 and variance is
    0.0342 (0.1849 std).}  This indicates that the variance is
  non-trivial and shows significant variation of the QoI with respect
  to the domain perturbation.
\end{remark}

\subsection{Sparse Grid convergence numerical experiment}
In this section the convergence rate of the sparse grid error is
tested without the truncation error. The purpose is to validate the
regularity of the solution with respect the stochastic parameters.

The mean $\mathbb{E}[Q]$ and variance $\var[Q]$ are computed with the
Clenshaw-Curtis isotropic {\it Sparse Grid Matlab Kit} \cite{Back2011}
for $N = 3,4,5$ dimensions.  The mean and variance are also computed
for $N = 3,4,5$ with a dimension adaptive sparse grid algorithm ({\it
  Sparse Grid Toolbox V5.1} \cite{Gerstner2003,spalg,spdoc}) and
Chebyshev-Gauss-Lobatto abscissas. In addition the following
parameters and experimental conditions are set:

\begin{itemize}
\item Let $a(\bmeta) = 1$ for all $\bmeta \in U$ and set the
  stochastic model parameters to $L = 19/50$, $L_P = 1$, $c= 1 /
  2.175$, $N = 15$,

\item The reference domain is discretized with a triangular mesh. The
  number of vertices are set in a $513 \times 513$ grid
  pattern. Recall that for the computation of the stochastic solution
  the fixed reference domain numerical method is used with the
  stochastic matrix $G(\by)$. Thus it is not necessary to re-mesh the
  domain for each perturbation.

\item The step size is set to $t_d := 1/1000$ and final time $T := 1$.

\item The QoI $Q(\hat u)$ is normalized by $Q(U)$ with respect the
  reference domain.

\end{itemize}

\begin{figure}[h]
  \begin{center}
  \begin{tabular}{cc}
\includegraphics[width=2.2in,height=2.2in]{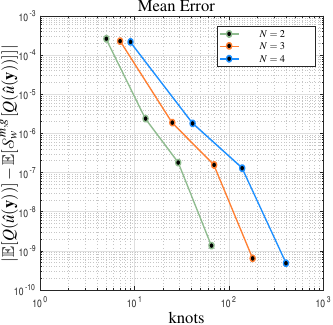} 
&
\includegraphics[width=2.2in,height=2.2in]{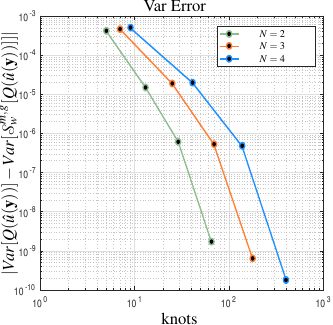} \\
\\
(a) & (b) \\
\end{tabular}
\end{center}
\caption{\cor{Isotropic sparse grid collocation convergence rates for $N =
  2,3,4$ with linear decay. (a) Mean error with respect to
  reference. Observe that the convergence rate is faster than
  polynomial.  (b) Variance error \cor{with respect to the
    reference}. For this case observe that the convergence rate is
  subexponential.}}
\label{results:fig2}
\end{figure}

In Figure \ref{results:fig2} (a) and (b) the normalized mean and
variance errors are shown for $N_s =2,3,4$. Each black marker
corresponds to a sparse grid level up to $w = 4$. For (a) we observe a
faster than polynomial convergence rate. Theoretically, the predicted
convergence rate should approach sub-exponential. This is not quite
clear from the graph as a higher level ($w \geq 5$) is needed to
confirm the results. However, this places the simulation beyond the
computational capabilities of the available hardware. In contrast, for
(b), the variance error convergence rate is clearly sub-exponential,
as the theory predicts.

\cor{
\begin{remark}
  In this work for simplicity we only demonstrate the application of
  isotropic sparse grids to the random domain problem. However,
  a significant improvement in error rates can be achieved by using an
  anisotropic sparse grid. By adapting the number of knots across
  each dimension to the decay rate of $\lambda_n$, $n =0,1,\dots,N$ a
  higher convergence rate can be achieved. In particular, if the decay
  rate of $\lambda_n$ is relative fast it will be not necessary to
  represent all the dimensions of $\Gamma$ to high accuracy.
  \end{remark}
}

\subsection{\cor{Truncation experiment}}

\begin{figure}[h]
\begin{center}
  \begin{tabular}{cc}
\includegraphics[width=2.2in,height=2.1in]{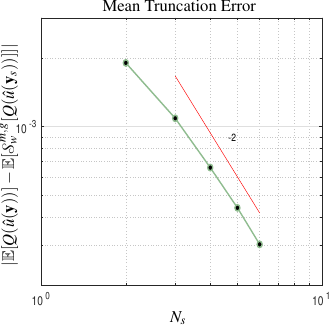}
&
\includegraphics[width=2.2in,height=2.1in]{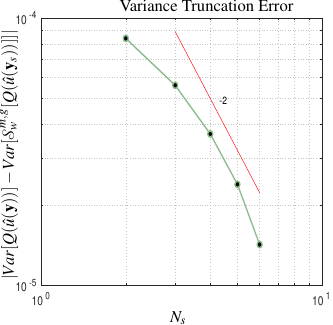}
\\
\corg{(a)} & \corg{(b)} \\
\end{tabular}
\end{center}
\caption{\corg{Truncation error with respect to the number of
    dimensions $N_s$ for linear decay $k = 1$. (a) Mean
    error. (b) Variance error. In both cases the decay appears
    quadratic, which is faster than the predicted convergence rate.}}
\label{results:fig3}
\end{figure}

\cor{The truncation error with respect to $N_s$ is analyzed and
  compared with respect to $Q(\hat{u}(\by))$ for $N = 15$ dimensions,
  $k = 1$ and $k = 1/2$.} \corg{The coefficient $c$ is changed to
  $1/4.35$.}  In Figure \ref{results:fig3} the truncation error is
plotted for (a) the mean and (b) the variance with respect to the
number of truncated dimensions $N_s$ for the linear decay $k = 1$.
From these plots observe that the convergence \cor{rates are close} to
quadratic, which is at least one order of magnitude higher than the
derived truncation convergence rate. \cor{Furthermore, in Figure
  \ref{results:fig4} the mean and variance error are shown for $k =
  1/2$.  As observed, the decay rate appears at least linear, which is
  at least twice the decay rate of the theoretical convergence rate.}
The numerical results \cor{shows that} in practice a higher
convergence rate is achieved than what the theory predicts.

\cor{
  \begin{figure}[h]
\begin{center}
  \begin{tabular}{cc}
\includegraphics[width=2.2in,height=2.1in]{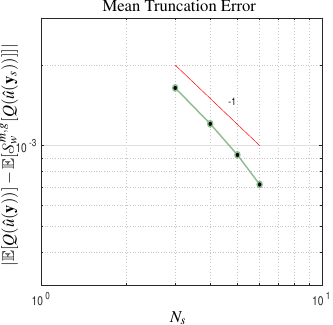}
&
\includegraphics[width=2.2in,height=2.1in]{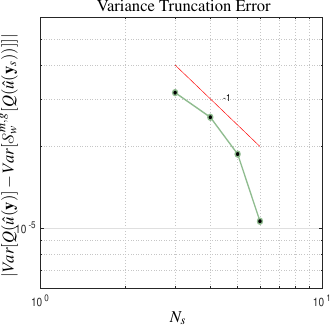}
\\
\corg{(a)} & \corg{(b)} \\
\end{tabular}
\end{center}
\caption{\corg{Truncation Error with respect to the number of dimensions
    $N_s$ for sqrt decay $k = 1/2$. (a) Mean error. (b) Variance
    error. In both cases the decay appears at least linear.}}
\label{results:fig4}
\end{figure}
}

\subsection{Forcing function truncation experiment}
\begin{figure}[h]
\begin{center}
  \begin{tabular}{cc}
\includegraphics[width=2.2in,height=2.1in]{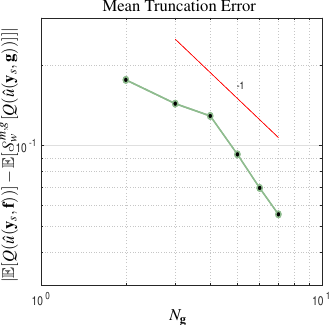}
&
\includegraphics[width=2.2in,height=2.1in]{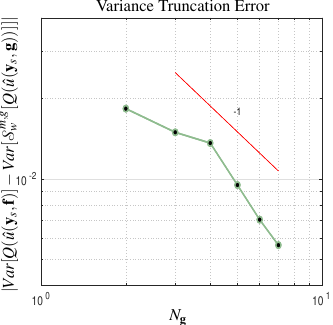}
\\
\corg{(a)} & \corg{(b)} \\
\end{tabular}
\end{center}
\caption{\corg{Forcing function truncation error with respect to the
    number of dimensions $N_{\bg}$. The decay of the coefficients
    $c(t,f_n)$, for $n = 1, \dots, N_{\bfq}$ are set to $1/n$.  The
    decay of the (a) Mean truncation error and the (b) Variance
    truncation error appears to be faster than linear.}}
\label{results:fig5}
\end{figure}

For the last numerical experiment the decay of the forcing function
truncation error (II) is tested with respect to the number of
dimensions $N_{\bg}$.  We compare the mean and variance error of
$Q(\bg,\by_s)$ with respect to $Q(\bfq,\by_s)$, where
\corg{
\[
f(\bx,\bfq,\by_s,t) = \sum_{n = 1}^{N_{\bfq}} c_{n}(t,f_n)
\xi_{n}(\bx,\by_s),\,\mbox{\&}\,
f(\bx,\bg,\by_s,t) = \sum_{n = 1}^{N_{\bg}} c_{n}(t,f_n)
\xi_{n}(\bx,\by_s),
\]
$\bx \in \Dw$ and $N_{\bfq} > N_{\bg}$. The maps $\xi_n:\Dw
\rightarrow 1$, for $n = 1, \dots, N$, are defined as
\[
\xi_n(x_1,x_2) := \exp\left(\frac{-(x_1 - a_n)^2}{\sigma}\right)
\exp\left(\frac{-(x_2 - b_n)^2}{\sigma}\right),
\]
where $\sigma = 0.001$.  The coefficients $a_n, b_n \in \R$ are given
such that $\xi_n$ are centered in a 4 by 4 grid.
Let $\ba := [ \frac{1}{4}\, \frac{5}{12},\, \frac{7}{12},\,
  \frac{3}{4}]$ $\bb := [ \frac{5}{8}\, \frac{17}{24},\,
  \frac{19}{24},\, \frac{7}{8}]$, then for $i = 1,\dots,4$ and $j =
1,\dots,4$ let $a_{4*(i-1)+j} := \ba[i]$, $b_{4*(i-1) + j} := \bb[j]$.
Furthermore,}
\begin{itemize}
\item For $n = 1,\dots,N_{\bfq}$, $f_{n}$ are
  independent uniform distributed in $(-\sqrt{3},\sqrt{3})$, and
  $c_{n}(t,f_n) = f^{2}_{n}/n$ (linear decay of the coefficients).
\item The stochastic PDE is solved \corg{on the domain $\Dw$} with a
  $513 \times 513$ triangular mesh.
\item \corg{$N_{\bfq} = 12$, $N_{s} = 2$, $N_{\bg} = 2, \dots, 7$ and
  $c = 1/4.35$.}
\item $\eset{Q(\by_s,\bfq)}$ and $var[Q(\by_s,\bfq)]$ are computed
  with a dimensional adaptive sparse grid with $\approx$ \corg{15,000}
  collocation points and a Chebyshev abscissa \cite{Gerstner2003}.
\item $\eset{Q(\by_s,\bg)}$ and $var[Q(\by_s,\bg)]$ are computed
  with the Clenshaw-Curtis isotropic {\it Sparse Grid Matlab Kit}
  \cite{Back2011} for $N_{\bg} = 2,\dots,7$.
\end{itemize}
By setting the coefficients to $c_{n}(t,f_n) = f^{2}_{n}/n$ we have a
non-linear mapping from the forcing function to the solution. From
Theorem \ref{erroranalysis:theorem2} the errors $|\mathbb{E}[Q(\hat
  u(\by_s,\bfq))] - \mathbb{E}[\mcS^{m,g}_w[Q(\hat u(\by_s,\bg))]]| $
and $|Var[Q(\hat u(\by_s,$ $\bfq)] - Var[\mcS^{m,g}_w[Q(\hat
    u(\by_s,\bg))]]|$ decay as
\[
\left(\sum_{n=N_{\bg}+1}^{N_{\bfq} }
\eset{c^{2}_{n}(t,f_n)}\right)^{1/2}
\sim
\frac{1}{N_{\bg}}.
\]
In Figure \ref{results:fig5} the error of the mean and variance are
plotted with respect to the number of dimensions $N_{\bg}$. The error
decay appears to be faster than the theoretically derived rate of
$\sim 1/N_{\bg}$.

\section{Conclusions}

In this paper a rigorous convergence analysis is derived for a sparse
grid stochastic collocation method for the numerical solution of
parabolic PDEs with random domains. The following contributions are
achieved in this work:

\begin{itemize}

\item An analysis of the regularity of the solution with respect to
  the parameters describing the domain perturbation \cor{shows that}
  an analytic extension onto a well defined region $ \Theta_{\beta}
  \times \cor{\msF} \subset \C^{N+N_{\bfq}}$ exists.

\item \cor{Error estimates in the energy norm for the solution and the
  QoI are derived for sparse grids with Clenshaw Curtis abscissas.}
  The derived subexponential convergence rate of the sparse grid is
  consistent with numerical experiments.

\item A truncation error with respect to the number of random
  variables is derived. Numerical experiments show a faster convergence
  rate.


\end{itemize}

This approach is well suited for a moderate number of stochastic
variables, but becomes impractical for large problems with an
isotropic sparse grid.  However, the approach described in this paper
can be \cor{easily extended} to anisotropic sparse grids
\cor{\cite{Schillings2013, nobile2008b}}.  Moreover, new approaches
such as quasi-optimal sparse grids \cite{Nobile2016} are shown to have
exponential convergence.

\section*{Acknowledgments} \corg{I appreciate the excellent feedback, comments,
suggestions and time from the reviewers of this article.}

\bibliographystyle{plain}
\bibliography{citations}



\end{document}